\documentclass[a4paper,11pt]{article}
\title{\bf{Stability conditions and 
crepant small resolutions}
}
\date{}
\author{Yukinobu Toda}
\usepackage{amsmath,amssymb}
\newcommand{\aA}{\mathcal{A}}
\newcommand{\bB}{\mathcal{B}}

\newcommand{\dD}{\mathcal{D}}
\newcommand{\eE}{\mathcal{E}}
\newcommand{\fF}{\mathcal{F}}
\newcommand{\gG}{\mathcal{G}}

\newcommand{\lL}{\mathcal{L}}
\newcommand{\mM}{\mathcal{M}}

\newcommand{\oO}{\mathcal{O}}
\newcommand{\pP}{\mathcal{P}}

\newcommand{\sS}{\mathcal{S}}
\newcommand{\tT}{\mathcal{T}}
\newcommand{\uU}{\mathcal{U}}

\newcommand{\wW}{\mathcal{W}}

\newcommand{\zZ}{\mathcal{Z}}
\newcommand{\lr}{\longrightarrow}

\newcommand{\Supp}{\mathop{\rm Supp}\nolimits}
\newcommand{\Hom}{\mathop{\rm Hom}\nolimits}

\newcommand{\dotimes}{\stackrel{\textbf{L}}{\otimes}}
\newcommand{\dR}{\mathbf{R}}

\newcommand{\Pic}{\mathop{\rm Pic}\nolimits}

\newcommand{\id}{\textrm{id}}

\newcommand{\ch}{\mathop{\rm ch}\nolimits}

\newcommand{\Ext}{\mathop{\rm Ext}\nolimits}
\newcommand{\Spec}{\mathop{\rm Spec}\nolimits}

\newcommand{\Coh}{\mathop{\rm Coh}\nolimits}

\newcommand{\Aff}{\mathop{\rm Aff}\nolimits}

\newcommand{\cneq}{\mathrel{\raise.095ex\hbox{:}\mkern-4.2mu=}}
\newcommand{\eqcn}{\mathrel{=\mkern-4.5mu\raise.095ex\hbox{:}}}

\newcommand{\Aut}{\mathop{\rm Aut}\nolimits}

\newcommand{\Ex}{\mathop{\rm Ex}\nolimits}

\newcommand{\Stab}{\mathop{\rm Stab}\nolimits}

\newcommand{\oPPer}{\mathop{\rm ^{0}Per}\nolimits}
\newcommand{\iPPer}{\mathop{\rm ^{-1}Per}\nolimits}

\newcommand{\pPPer}{\mathop{\rm ^{\mathit{p}}Per}\nolimits}

\newcommand{\Imm}{\mathop{\rm Im}\nolimits}
\newcommand{\Ree}{\mathop{\rm Re}\nolimits}
\newcommand{\FM}{\mathop{\rm FM}\nolimits}

\newcommand{\Stabn}{\Stab _{\rm{n}}}

\newcommand{\length}{\mathop{\rm length}\nolimits}

\usepackage{array}
\usepackage{amscd}
\usepackage[all]{xy}

\newtheorem{thm}{Theorem}[section]
\newtheorem{prop}[thm]{Proposition}
\newtheorem{lem}[thm]{Lemma}
\newtheorem{defi}[thm]{Definition}
\newtheorem{rmk}[thm]{Remark}
\newtheorem{cor}[thm]{Corollary}
\newtheorem{step}{Step}
\newtheorem{sstep}{Step}

\newtheorem{prop-defi}[thm]{Proposition-Definition}

\begin{document}
\maketitle
\begin{abstract}
In this paper, we describe the spaces of stability conditions 
on the triangulated categories associated to   
three dimensional crepant small resolutions. 
The resulting spaces have chamber structures such that
each chamber   
corresponds to a birational model together with a 
special Fourier-Mukai transform.  
We observe that these spaces are covering spaces over certain 
open subsets of finite dimensional vector spaces, and 
determine their deck transformations. 
\end{abstract} 
\section{Introduction}

$\quad$

For a triangulated category $\dD$, the notion of stability 
conditions on $\dD$ was introduced by T.Bridgeland~\cite{Brs1}
 to give a mathematical 
framework for Douglas's notion of $\Pi$-stability~\cite{Dou1}, \cite{Dou2}. 
Roughly speaking, Bridgeland's stability condition $(Z, \pP)$
 on $\dD$
 consists of a group homomorphism $Z$ and full subcategories, 
 $$Z\colon K(\dD) \longrightarrow \mathbb{C}, \quad \pP(\phi)\subset 
 \dD \quad (\phi \in \mathbb{R}),$$
where $K(\dD)$ is the Grothendieck group of $\dD$, and 
this pair satisfies some axioms.
Given data as above, the objects of the subcategory $\pP (\phi)$ 
are called semistable of phase $\phi$, and
this gives a generalization of the classical notion of semistable
sheaves with a fixed slope on a smooth projective curve. 
In~\cite{Brs1}, T.Bridgeland showed that the set of locally finite 
stability conditions $\Stab (\dD)$ has a natural topology, and 
proposed to study this space as
a new categorical invariant. 
After he wrote the paper~\cite{Brs1}, some examples have been computed 
in~\cite{Brs2}, \cite{Brs3}, 
\cite{Brs4}, \cite{Tho}, \cite{Oka}, \cite{Mac}, \cite{IUU}, \cite{Ber}. 
The purpose of this paper is 
to describe the spaces of stability conditions 
on some triangulated categories associated 
to three dimensional crepant small resolutions.

\subsection{Crepant small resolutions}

Let
$Y=\Spec R$ for a noetherian local complete Gorenstein
$\mathbb{C}$-algebra $R$
of dimension three, and $0\in Y$ be the closed point. 
  Assume that $Y$ admits a resolution of singularities
 $$f\colon X\to Y=\Spec R$$ which is
an isomorphism in 
codimension one. 
Then $f\colon X\to Y$ is called a \textit{crepant small resolution}. 
It is well-known that $\omega _X=f^{\ast}\omega _Y$ and the 
exceptional locus $C=\cup _{i=1}^{N}C_i$
is a tree of 
rational curves.
Let $D^b (X)$ be the bounded derived category of coherent sheaves on $X$. 
We define the triangulated subcategory $\dD_{X/Y}\subset D^b (X)$ as follows:
$$
\dD_{X/Y} \cneq \{ E\in D^b (X) \mid \Supp (E)\subset C \}.$$
The purpose of this paper is to
 study the space of stability conditions on $\dD _{X/Y}$. 
Define the space $\Stab (X/Y)$ by 
$$\Stab (X/Y)\cneq \Stab (\dD _{X/Y}).$$
The simplest example is given by the morphism, 
$$f\colon X\longrightarrow Y=\Spec \mathbb{C}[[x,y,z,w]]/(xy-zw)$$
which is the blowing up by the ideal $(x,z)\subset \oO _Y$. 
In this case, the exceptional locus is a single rational curve, whose 
normal bundle is $\oO _C(-1)\oplus \oO _C(-1)$. Taking the blow-up 
by the ideal $(x,w)\subset \oO _Y$ gives another resolution 
$g \colon W \to Y$, and the diagram 
$$X \stackrel{f}{\longrightarrow} Y 
\stackrel{g}{\longleftarrow}W,$$
is an example of a $\textit{flop}$~\cite[Definition 6.10]{KM}.
A new feature in studying $\Stab (X/Y)$ is that
one has to take account of all the crepant small resolutions of 
$Y$, such as flops.
 We would like to explain this 
 in two different contexts, 
string theory and minimal model theory. 

\subsection{Viewpoint from string theory}
Here we give a rough picture on the space $\Stab (X/Y)$,
from the viewpoint of string theory. 
We use notions of string theory, 
for example stringy K$\ddot{\textrm{a}}$hler moduli spaces, and SCFT. 
These notions are explained in~\cite[Section 2]{Brs6}. 
For instance assume $X$ is a projective Calabi-Yau 3-fold, 
$\dD =D^b (X)$, 
and $\mM _K(X)$ is the stringy K$\ddot{\textrm{a}}$hler moduli 
space~\cite[(2.2)]{Brs6}. The space $\mM _K(X)$ is a subspace 
of the moduli space of SCFT, and the associated topological B-twists
 are unchanged 
along $\mM _K(X)$. Moreover if
 $\hat{X}$ is a mirror of $X$, $\mM _K(X)$ is supposed to
be isomorphic to the moduli space of complex structures on $\hat{X}$. 
 The space of stability conditions was introduced in order to understand 
 $\mM _K(X)$ mathematically. More precisely, it
  is believed that the quotient space of 
 $\Stab (\dD)$ by the actions of $\Aut (\dD)$ and $\mathbb{C}$
 contains $\mM _K(X)$~\cite[Remark 3.9]{Brs6}.  
 
 The moduli space of complexified K$\ddot{\textrm{a}}$hler forms 
 $\beta +i \omega \in H^2 (X, \mathbb{C})$ forms an open
  subset $\uU _X \subset \mM _K(X)$, which physicists call the
  \textit{neighborhood of 
  the large volume limit}. The important point is that
  there might be another topologically distinct Calabi-Yau 3-fold $W$
   such
  that $\uU _W$ is 
  contained in $\mM _K(X)$ in the moduli space of SCFT. 
   In this case, the associated B-twists 
  from the data of $X$ and $W$ are equivalent, therefore their 
  categories of $D$-branes are also equivalent. 
  Mathematically this means that there 
  exists an equivalence of the bounded derived categories of 
  coherent sheaves, 
  $$\Phi \colon D^b (W) \longrightarrow D^b (X).$$
  In the language of 
  algebraic geometry, $W$ is called a \textit{Fourier-Mukai partner}
  of $X$ and a flop gives one 
  example~\cite{B-O2},\cite{Br1},\cite{Ch},\cite{Ka1}.
  Therefore in describing $\mM _K(X)$, we have to take into account
  the neighborhoods of the large volume limits corresponding to
   several Fourier-Mukai 
  partners. 
  
  In fact in
   the case of a flop $X \to Y \leftarrow W$, 
  P.Aspinwall~\cite[Figure 2]{Asp}
   describes the localized picture of $\mM _K(X)$
  assuming all the curves in $X$ except the flopping curve $C\subset X$
  are quite ``big". The resulting picture is a 2-sphere minus three points, 
  $\uU _X$, $\uU _W$ are disjoint and their union is dense in the sphere. 
  The string theory has 
  singularities at one of the deleted points, and the other two points 
  are large volume limits corresponding 
  to $X$ and $W$. 
  Thus one should have a similar picture in the context of 
  Bridgeland's stability conditions.
  It seems that
  our localized category $\dD _{X/Y}$ and the space $\Stab (X/Y)$ 
  provide the right framework, 
   thus we expect $\Stab (X/Y)$ is described 
   via 
  Fourier-Mukai partners such as flops.

\subsection{Viewpoint from minimal model theory}

Our interest also comes from the
birational geometry, especially minimal model program, 
simply MMP~\cite{KM},~\cite{KMM}. 
The MMP is a program aimed to find a good birational model for 
a given projective variety, by contracting extraneous rational curves. 
The output is either a minimal model or Mori fiber space. 
 
One of the points where the 
 three dimensional MMP differs from the two dimensional one 
is that birational 
minimal models are not necessary unique,  
but connected by a sequence of flops~\cite[Theorem 6.38]{KM}. 
The philosophy of
 Y.Kawamata~\cite{Ka2} is that one can capture the set of 
 birational minimal models via a chamber structures on 
 the \textit{movable cone}~\cite[Definition 1.1]{Ka2}. 
 According to~\cite[Theorem 2.3]{Ka2}, chambers on 
 the movable cone are given by the ample cones 
 of birational minimal models. 
 
 Now let us consider two birational
 three dimensional minimal models $W
\dashrightarrow X$. Then there exists an 
equivalence of bounded 
derived categories
of coherent sheaves~\cite{B-O2},\cite{Br1},\cite{Ch},\cite{Ka1}, 
$$\Phi \colon D^b (W) \longrightarrow D^b (X).$$
Let $\Stab (X)$ be the space of stability conditions on $D^b (X)$. 
As a substitute for the ample cone, one expects to find
 a certain region 
$U_X \subset \Stab (X)$, which corresponds to the neighborhood of 
the large volume limit in string theory. 
Then one can transfer stability conditions in $U_{W}\subset \Stab (W)$
by the equivalence $\Phi$ to get 
the region $U(W, \Phi)\subset \Stab (X)$. 
 In summary for a birational minimal model 
$W$ and an equivalence $\Phi$ as above, one obtains the
correspondence, 
$$(W, \Phi) \longmapsto U(W, \Phi)
 \subset \Stab (X).$$ 
This picture is quite similar to the picture of the 
movable cone~\cite[Theorem 2.3]{Ka2}. 
Thus we
guess the existence of the chamber structure on $\Stab (X)$,  
which enables us to capture the pair $(W, \Phi)$ as above. 
  Unfortunately, there are some technical 
issues in working with
Bridgeland's stability conditions for the derived categories of 
projective 3-folds. In particular we do not know how to
construct examples of stability conditions in this case.  (Also see
the last part of~\cite[Section 4]{Brs6}.)
Despite this problem, 
our category $\dD _{X/Y}$ is quite amenable to
studying stability conditions, and also 
sufficient for realizing our purpose.

\subsection{The main results}

Let $f\colon X\to Y$ be a crepant small resolution as in the first part 
of this introduction.   
For some technical reasons, we put the following additional 
assumption, 
\begin{itemize}
\item There is a hyperplane section $0\in Y_0 \subset Y$ such that its 
pull-back $X_0 \cneq f^{-1}(Y_0)$ is smooth.
\end{itemize}

First we give the standard region in Definition~\ref{LV} below, 
$$U_X \subset \Stab (X/Y),$$
which should correspond to the neighborhood of the large volume limit. 
Unfortunately, $U_X$ is not an open subset of $\Stab (X/Y)$, but
 open in the subspace of $\Stab (X/Y)$ which 
we call the \textit{normalized stability conditions}, 
$$\Stabn (X/Y) \cneq \{ \sigma =(Z, \pP) \in \Stab (X/Y) \mid 
Z([\oO _x])=-1 \}.$$
In this paper, we give a description of the normalized version of
the stability conditions. Because other stability conditions are 
obtained by the actions of the additive group $\mathbb{C}$
 from the normalized
stability conditions, it is enough to study $\Stabn (X/Y)$ for
our purpose. Here the action of $\mathbb{C}$ is as follows:
for $\lambda \in \mathbb{C}$ and $\sigma =(Z, \pP)\in \Stab (X/Y)$, 
$\lambda (\sigma)=(Z',\pP ')$ with $\pP ' (\phi )=\pP (\phi +\Ree \lambda)$
and $Z' (E)=\exp (-i\pi \lambda)Z(E)$.
Let $\Stabn ^{\circ}(X/Y)$ be the connected component of $\Stabn (X/Y)$ 
which contains $U_X$. 

Next we introduce the set $\FM (X)$ to be the set of 
pairs $(W, \Phi)$, where $g\colon W \to Y$ is a crepant 
small resolution, and 
$$\Phi \colon D^b (W) \longrightarrow D^b (X)$$
gives an equivalence of derived categories, such that 
$\Phi$ is given in a special way, as described in Definition~\ref{FM} below. 
By transferring the open set $U_{W}\subset \Stabn (W/Y)$ through the 
equivalence $\Phi$, we give the following open set in Definition~\ref{other},
$$U_{}(W, \Phi) \subset \Stabn (X/Y).$$

The following theorem realizes our purpose, that is
 the existence of the desired chamber structure on
the space of stability conditions.

\begin{thm}\label{cham}
We have the following union of chambers:
$$\mM _{}\cneq \bigcup _{(W,\Phi)\in \rm{FM}(X)}
U_{}(W, \Phi _{})\subset \Stabn^{\circ}(X/Y), $$
such that two chambers are either disjoint or equal. Moreover
we have $\overline{\mM}_{}=\Stabn^{\circ}(X/Y)$. 
\end{thm}

By using the chamber structure given in Theorem~\ref{cham}, 
we show $\Stabn ^{\circ}(X/Y)$ is a regular covering space over a
certain open subset of a finite dimensional complex vector space.
Let $\Lambda _f$ be the root lattice associated to the 
exceptional locus of $f$. 
We denote by $N^1 (X/Y)$ the group of numerical classes of 
$\mathbb{R}$-divisors and $N^1 (X/Y)_{\mathbb{C}}\cneq 
N^1 (X/Y)\otimes _{\mathbb{R}}\mathbb{C}$.
 Then the elements of $N^1 (X/Y)_{\mathbb{C}}$ are regarded
as complex functions on $\Lambda _f$. 
Let $V(\Lambda _f)$ be the set of roots of $\Lambda _f$. 
For a 
root $v\in V(\Lambda _f)$ and $k\in \mathbb{Z}$,
 define $H_v$ as follows:
 $$H_v \cneq \{ \beta +i \omega \in N^1 (X/Y)_{\mathbb{C}} \mid
 (\beta +i \omega) v \in \mathbb{Z} \}. $$

\begin{thm}\label{cham2}
The map
$$\zZ _{X} 
\colon \Stabn ^{\circ}
(X/Y) \longrightarrow N^1 (X/Y)_{\mathbb{C}}\setminus 
\bigcup _{v\in V(\Lambda _f)}H_v , $$
defined by sending a stability condition to its central charge, 
is a regular covering map. 
\end{thm}

The results of this paper together with some developments on 
derived categories and stability conditions 
are reviewed in T.Bridgeland's manuscript~\cite{ICM} for 
the ICM talk in 2006.

\subsection*{Acknowledgement}
This paper was written while the author was visiting 
the University of Sheffield from September 2005 to 
December 2005.  
The author thanks Tom Bridgeland for leading the author
to the stability conditions, and many useful discussions and 
comments. 
Especially the idea of constructing the stability conditions in 
Lemma~\ref{nor} is due to him. 
The author also thanks his advisor professor Yujiro Kawamata, who 
recommended him to visit the University of Sheffield.  
Finally the author thanks the referee for reading the manuscript 
carefully, and giving him several nice comments to make the manuscript 
much readable. 
The author
 is supported by the Japan Society for the Promotion of Sciences Research 
Fellowships for Young Scientists, No 1611452.

\section{Generalities}
In this paper, all the schemes are defined over $\mathbb{C}$. 
For a scheme $X$, we denote 
by $\Coh (X)$ and $D^b (X)$ the Abelian category of 
coherent sheaves and its bounded derived category respectively. 
The shift functor on $D^b (X)$ is denoted by $[1]$. 
Also for a subscheme $Z\subset X$, we denote by 
 $\Coh _Z (X) \subset \Coh (X)$ the 
 subcategory whose objects 
 are supported on $Z$.
 For an object $E\in D^b (X)$, its support is defined by 
 $$\Supp (E)\cneq \bigcup _{p\in \mathbb{Z}} \Supp H^p(E) \subset X.$$
For a triangulated category $\dD$, its $K$-group is denoted by 
$K(\dD)$. 
\subsection{Stability conditions on triangulated categories}
$\quad$
The notion of stability conditions on triangulated categories 
is introduced in~\cite{Brs1} to give the mathematical 
framework for the Douglas's work on $\Pi$-stability~\cite{Dou1}, \cite{Dou2}. 
Here we collect the basic definitions and results in~\cite{Brs1}.

\begin{defi}
A stability condition on a triangulated category $\dD$ 
consists of data $\sigma =(Z, \pP)$, 
where $Z\colon K(\dD)\to \mathbb{C}$ is a linear map, 
and $\pP (\phi)\subset \dD$ is a full additive subcategory 
for each $\phi \in \mathbb{R}$, 
which satisfy the following:

 \begin{itemize}
 \item $\pP (\phi +1)=\pP (\phi)[1].$ 
\item  If $\phi _1 >\phi _2$ and $A_i \in \pP (\phi _i)$, then 
$\Hom (A_1, A_2)=0$. 
\item  If $E\in \pP (\phi)$ is non-zero,
 then $Z(E)=m(E)\exp (i\pi \phi)$ for some 
$m(E)\in \mathbb{R}_{>0}$. 
\item For a non-zero object $E\in \dD$, we have the 
following collection of triangles:
$$\xymatrix{
0=E_0 \ar[rr]  & &E_1 \ar[dl] \ar[rr] & & E_2 \ar[r]\ar[dl] & \cdots \ar[rr] & & E_n =E \ar[dl]\\
&  A_1 \ar[ul]^{[1]} & & A_2 \ar[ul]^{[1]}& & & A_n \ar[ul]^{[1]}&
}$$
such that $A_j \in \pP (\phi _j)$ with $\phi _1 > \phi _2 > \cdots >\phi _n$. 
\end{itemize}
\end{defi} 
Here $Z$ is called the \textit{central charge}. 
Each $\pP (\phi)$ is an Abelian category,
the non-zero objects of $\pP (\phi)$ are called \textit{semistable of 
phase} $\phi$, and simple objects of $\pP (\phi)$ are called 
\textit{stable}. The objects $A_j$ are called \textit{semistable factors}
 of $E$ with 
respect to $\sigma$. 
The following proposition is useful in constructing stability conditions. 
\begin{prop}\emph{\bf{\cite[Proposition 4.2]{Brs1}}}\label{tstru}
Giving a stability condition on $\dD$ is equivalent to giving a heart of a
bounded t-structure $\aA \subset \dD$, and a group homomorphism 
$Z\colon K(\dD)\to \mathbb{C}$ called the stability function,
 such that for a non-zero object $E\in \aA$ one has
$$ Z(E)\in
\{ r\exp (i\pi \phi) \mid r>0, 0<\phi \le 1\},$$
and the pair $(Z, \aA)$ 
satisfies the Harder-Narasimhan 
property.
\end{prop}
For the Harder-Narasimhan property, we refer~\cite[Definition 2.3]{Brs1}. 
Also a necessary condition
 for this property are found in~\cite[Proposition 2.4]{Brs1}. 
Given data $(Z, \aA)$ as above and a non-zero object 
$E\in \aA$, we define $\phi (E)\in (0,1]$ such that 
$Z(E)\in \mathbb{R}_{>0}\exp (i\pi \phi (E))$ uniquely. 
We call $\phi (E)$ the \textit{phase} of $E$.
\begin{rmk}\label{filt}\emph{
The condition (b) of~\cite[Proposition 2.4]{Brs1}
is satisfied if $\aA$ is noetherian, i.e. 
for $E\in \aA$, there are no infinite sequences of subobjects in $\aA$, 
$$E_1 \subset E_2 \subset \cdots E_j \subset \cdots \subset E.$$
Therefore in order to check the Harder-Narasimhan property, 
it is enough to check that $\aA$ is noetherian and the 
condition (a) of~\cite[Proposition 2.4]{Brs1}, i.e. there 
are no infinite sequence of subobjects in in $\aA$, 
$$\cdots \subset E_{j+1}\subset E_j \subset \cdots \subset E_2 \subset E_1,$$
with $\phi (E_{j+1})>\phi(E_j)$ for all $j$. }
\end{rmk}

The set of stability conditions which satisfies the technical condition
\textit{local finiteness}~\cite[Definition 5.7]{Brs1}
is denoted by $\Stab (\dD)$. 
It is shown in~\cite[Section 6]{Brs1} that $\Stab (\dD)$ has a 
natural topology. 
Forgetting the information of $\pP$, we have the map
$$\zZ \colon \Stab (\dD) \longrightarrow \Hom _{\mathbb{Z}}(K(\dD), 
\mathbb{C}).$$
\begin{thm}\emph{\bf{\cite[Theorem 1.2]{Brs1}}}\label{lois}
For each connected component $\Sigma \subset \Stab (\dD)$, 
there exists a linear subspace $V(\Sigma )\subset \Hom _{\mathbb{Z}}(K(\dD), 
\mathbb{C})$, such that $\zZ$ restricts to a local homeomorphism, 
$\zZ \colon \Sigma \to V(\Sigma)$. 
\end{thm}
In general $\Stab (\dD)$ is infinite dimensional, 
so we usually consider only numerical stability conditions 
as in~\cite{Brs1},~\cite{Brs2}. 
 But if we assume 
$K(\dD)$ is finitely generated,
Theorem~\ref{lois} implies each connected component
$\Sigma \subset \Stab (\dD)$ is a complex manifold.

\section{Geometry on crepant small resolutions}
Let $f\colon X\to Y=\Spec R$ be a three dimensional 
crepant small resolution as in the introduction. 
The exceptional locus $C\subset X$ is a tree of rational 
curves, 
$$C=C_1 \cup C_2 \cup \cdots \cup C_N,$$
with each $C_i$ isomorphic to $\mathbb{P}^1$. 
 (See for example~\cite[Lemma 3.4.1]{MVB}.) 
In this paper, we put the 
following 
additional assumption. 
\begin{itemize}
\item 
There exists a hyperplane section $0\in Y_0 \subset Y$
such that 
its pull-back $X_0 \cneq f^{-1}(Y_0)$ is smooth. 
\end{itemize}
In general, it is known that a general hyper plane section 
$Y_0$ is a rational double point, 
and $X_0 \to Y_0$ is a partial resolution, i.e. the
minimal resolution $X_0 ' \to Y_0$ factors through 
$X_0 ' \to X_0 \to Y_0$. (See~\cite[(1.1), (1.14)]{Rei}.) 
Also it is known that for a given crepant small resolution 
$f\colon X\to Y$, there exists a finite map $Y' \to Y$ such that 
$Y'$ admits a crepant small resolution which satisfies the above 
assumption. (See~\cite[Theorem 4.28]{KM}).
The above assumption will be required 
in Lemma~\ref{pers}, Subsection~\ref{root} and Lemma~\ref{ch}.
Of course this assumption is satisfied 
in the case of the resolution of the ordinary double point described in the 
introduction. 
As in the introduction, we define $\dD _{X/Y}$ to be 
$$\dD _{X/Y}\cneq \{ E\in D^b (X) \mid \Supp E \subset C \}.$$
We collect some notation and known results on this resolution, 
and gives some lemmas. All the lemmas in this section will be 
proved in Section 6. 

\subsection{Terminology and results from birational geometry}
Here we introduce standard terminology in 
birational geometry, for example used in~\cite[Definition 1.1]{Ka2}. 
Two divisors $D_1$, $D_2$ on $X$ are called \textit{numerically equivalent} 
over $Y$ if and only if 
$D_1 \cdot C_i =D_2 \cdot C_i$ for all $1\le i\le N$.
Similarly, two 
 1-cycles $Z_1$, $Z_2$ on $X$ contracted by $f$
 are \textit{numerically equivalent} 
if and only if $D\cdot Z_1=D\cdot Z_2$ for every divisor $D$ on $X$. 
\begin{defi}
We define the $\mathbb{R}$-vector spaces $N^1 (X/S)$, $N_1(X/S)$ to be
\begin{align*}
N^1 (X/S) & \cneq \{ \emph{Divisors on }X \} / (
\emph{numerical equivalence over }Y ) \otimes _{\mathbb{Z}}
 \mathbb{R}, \\
 N_1 (X/S) & \cneq \{ \emph{One cycles on }X
  \emph{ contracted by }f \} / (
\emph{numerical equivalence}) \otimes _{\mathbb{Z}}\mathbb{R}
\end{align*}
\end{defi}
By the definition, one has the perfect pairing, 
$$N^1 (X/Y) \times N_1 (X/S) \ni (D, Z) \longmapsto D\cdot Z \in \mathbb{R}.$$
Moreover since we are assuming 
$Y$ is complete, there exist divisors $D_i$ on $X$
for $1\le i\le N$ such that 
$$D_i \cdot C_j = \left\{ \begin{array}{ll} 1 \quad (i=j) \\ 0 \quad 
(i\neq j) 
\end{array} \right. $$ 
by~\cite[Lemma 3.4.4]{MVB}. Therefore we have
$$N^1 (X/S)=\bigoplus _{1\le i\le N} \mathbb{R}[D_i], \quad 
N_1 (X/S)=\bigoplus _{1\le i \le N}\mathbb{R}[C_i],$$
in our case. Similarly we introduce the one dimensional $\mathbb{R}$-vector 
spaces, 
$$N^0 (X/Y)=\mathbb{R}[X], \quad N_0 (X/Y)=\mathbb{R}[p], $$
for a closed point $p\in C$ as \textit{numerical classes of codimension 
zero cycles} and \textit{zero dimensional cycles}. We have the pairing, 
$$ N^0 (X/Y) \times N_0 (X/Y) \ni (a[X], b[p]) \longmapsto 
ab \in \mathbb{R}.$$
Let $N^1 (X/Y)_{\mathbb{C}}\cneq N^1 (X/Y)\otimes _{\mathbb{R}}\mathbb{C}$.
\begin{defi}\label{amp}
We define the ample cone $A(X/Y)$ and the
complexified ample cone $A(X/Y)_{\mathbb{C}}$ to be
\begin{align*} A(X/Y) & \cneq \{ \emph{Numerical classes of ample }\mathbb{R} 
\emph{-divisors } \} \subset N^1 (X/Y) \\
& = \bigoplus _{1\le i\le N} \mathbb{R}_{>0}[D_i], \\
A(X/Y)_{\mathbb{C}} & \cneq \{ \beta +i \omega \in N^1 (X/Y)_{\mathbb{C}} 
\mid \omega \in A(X/Y) \}. \end{align*}
Also for each $1\le j \le N$ and $k\in \mathbb{Z}$,
 define $\wW _j$ and $\wW _{j,k}$ to be 
\begin{align*}
\wW _j & \cneq \{ \beta +i \omega \in N^1 (X/Y)_{\mathbb{C}} \mid 
\omega \cdot C_j =0 \emph{ and }\omega \cdot C_{j'} >0
 \emph{ for } j\neq j' \}\\
\wW _{j,k}&\cneq \{ \beta +i \omega \in \wW _j \mid \beta \cdot C_j \in 
(k-1, k) \}. \end{align*}
\end{defi}
 Note that the union $\bigcup _{1\le j\le N}\wW _j$ is 
 the codimension one boundary of the \textit{complexified nef cone}
  $\overline{A}(X/Y)_{\mathbb{C}}
 \subset N^1 (X/Y)_{\mathbb{C}}$. 
 For each $1\le i \le N$, some multiple of the
  divisor 
 $$\sum_{j\neq i}D_j,$$
 is base point free by~\cite[Theorem 3.3]{KM}, thus there exists a 
 birational contraction $g_i \colon X \to Y_i$ which 
 contracts only $C_i$. Then $\wW _i$ is written as 
 $$\wW _i =g_i^{\ast}\aA (Y_i /Y)_{\mathbb{C}}.$$
 Furthermore one can construct its \textit{flop}~\cite[Theorem 6.14]{KM}
 and obtain the diagram below:
\begin{align} \xymatrix{
(C_i \subset X) \ar[r]^{g_i} \ar[dr]_{f} & (p_i \in Y_i )\ar[d]^{h_i}
 & (X_i ^{\dag} \supset C_i ^{\dag})\ar[l]_{g_i ^{\dag}} \ar[ld]^{f_i ^{\dag}} \\
 & Y. & }\end{align}
 Next let $g\colon W\to Y$ be another crepant small resolution, 
 and $\phi =f^{-1}\circ g \colon W\dashrightarrow X$ be the birational 
 map. Because $\phi$ is an isomorphism in codimension one, one
 has the isomorphism of the groups of divisor classes, 
  $\phi _{\ast}\colon N^1(W/Y) 
 \to N^1 (X/Y)$ called the \textit{strict transform}. 
 We use the same notation $\phi _{\ast}\colon \Pic (W) \to \Pic (X)$ 
 for the isomorphism of Picard groups. 
 The following theorem gives a chamber structure on $N^1 (X/Y)$, 
 where each chamber corresponds to the ample cone of a
  crepant small resolution
 of $Y$. 
 
 \begin{thm}\emph{\bf{\cite[Main theorem]{KaMa}}},
  \emph{\bf{\cite[Theorem 2.3]{Ka2}}}
 \label{model}
 The number of crepant small resolutions $g\colon W\to Y$ is finite 
 up to isomorphism. One has the decomposition, 
 $$N^1 (X/Y)=\bigcup _{(W, \phi)}\phi _{\ast}\overline{A}(W/Y).$$
 Here $(W, \phi)$ is a pair of a crepant small resolution
 $g\colon W\to Y$ and the birational map $\phi =f^{-1}\circ g$. 
 Moreover $\phi _{\ast}A(W/Y)\cap \phi '_{\ast}A(W'/Y)\neq \emptyset$
 if and only if there exists an isomorphism $h\colon W' \to W$ such 
 that $\phi \circ h =\phi '$. 
 \end{thm}

\begin{rmk}\emph{In our case, $f\colon X\to Y$ is an isomorphism in codimension one. 
Therefore the $f$-effective $f$-movable cone and $f$-effective $f$-nef 
cone defined in~\cite[Definition 1.1]{Ka2} coincide with
$N^1 (X/Y)$ and $\overline{A}(X/Y)$ respectively.}
\end{rmk}

\subsection{Chern characters and Riemann-Roch theorem}
We consider the chern characters which take values in our 
vector spaces $N^{\ast}(X/Y)$ and $N_{\ast}(X/Y)$, 
\begin{align*}
D^b (X) \ni E & \longmapsto (\ch _0 (E), \ch _1 (E)) \in 
N^0 (X/Y) \oplus N^1 (X/Y), \\
\dD _{X/Y} \ni F & \longmapsto (\ch _2 (F), \ch _3(F) ) \in 
N_1 (X/Y) \oplus N_0 (X/Y).
\end{align*}
Usually, chern characters take values in Chow groups. But since 
rationally equivalent two cycles are numerically equivalent, 
$\ch _0$, $\ch _1$
can take values in $N^0 (X/Y)$, $N^1 (X/Y)$. 
Also $F\in \dD _{X/Y}$ is supported on $C$, thus
$\ch _2$, $\ch _1$ take values in $N_1 (X/Y)$, $N_0 (X/Y)$. 

Let us take $E\in D^b (X)$, $F\in \dD _{X/Y}$. Since 
the support of $F$ is proper, the space $\Ext _{X}^i (E, F)$ 
is finite dimensional, and zero except for a finite number of 
$i\in \mathbb{Z}$. Therefore the number
$$\chi (E, F)\cneq \sum _{i\in \mathbb{Z}}(-1)^i \dim \Ext ^i _X 
(E, F) \in \mathbb{Z}, $$
makes sense. The Riemann-Roch theorem~\cite[Corollary 18.3.1]{Fu}
 implies,
$$\chi (E, F) =\ch _0 (E)\cdot \ch _3(F) -\ch _1(E)\cdot \ch _2(F).$$

\subsection{Perverse t-structures}
Here we introduce the Abelian categories $\pPPer (X/Y)\subset D^b (X)$ 
for $p=0, -1$, which 
are introduced in~\cite{Br1} to construct the flops. 
These categories provide a central technical tool in this paper. 
For the precise definition, we refer~\cite[Section 3]{Br1} and~\cite[Section 3]{MVB}. 
First we give the known results of~\cite{Br1},\cite{MVB}.
Let $S_0$, $S_0'$, and $S_i$ for $1\le i\le N$ be
 $$S_0 \cneq \omega _C [1], \quad S_0 ' \cneq \oO _C, \quad 
 S_i \cneq \oO _{C_i}(-1).$$
 Here the scheme structure on $C_i$ is reduced in the 
 definition of $S_i$ for $1\le i\le N$. However 
 the scheme structure on $C$ in the definition of 
 $S_0$, $S_0'$ must be the scheme theoretic fiber of $f\colon X\to Y$, 
 thus possibly non-reduced. 

\begin{lem}\emph{\bf{\cite[Lemma 3.2]{Br1}}}, 
\emph{\bf{\cite[Section 3]{MVB}}}
\label{per}
For $p=0$ or $-1$, there exist hearts of bounded t-structures 
$\pPPer (X/Y) \subset D^b (X)$, such that 
$E\in D^b (X)$ belongs to $\pPPer (X/Y)$
if and only if the following holds: 

(i) $E$ is concentrated in degrees $[-1, 0]$.  

(ii) $f_{\ast}H^{-1}(E)=0$ and $R^1 f_{\ast}H^0 (E)=0$. 

(iii) For $j\ge 1$, we have the following:
$$\Hom (S_j, H^{-1}(E))=0, \quad (p=0), \qquad 
\Hom (H^0 (E), S_j)=0, \quad (p=-1).$$

\end{lem}

Here we define the categories $\pPPer (\dD _{X/Y})$ for $p=-1, 0$ 
to be 
$$\pPPer (\dD _{X/Y})\cneq \pPPer (X/Y) \cap \dD _{X/Y}.$$
M.Van den Bergh~\cite{MVB} determined the simple objects of these categories.

\begin{prop}\emph{\bf{\cite[Proposition 3.5.8]{MVB}}}\label{simple}
The categories 
$$\oPPer (\dD_{X/Y}),  \quad \mbox{ and } \quad
\iPPer (\dD _{X/Y})$$ are hearts of bounded t-structures on $\dD _{X/Y}$, 
and are finite length 
Abelian categories. Their simple objects are 
$$\{ S_i \}_{0\le i \le N},  \mbox{ and }
\quad \{ S_0 ', S_i [1]\} _{1\le i\le N}$$ respectively. 
 \end{prop} 

Let $g_i \colon X \to Y_i$ be the contraction which 
contracts $C_i$, as in the diagram (1). Note that one 
can define $\pPPer (X/Y_i)\subset D^b (X)$, and 
$\pPPer (\dD _{X/Y_i})\subset \dD _{X/Y_i}$ similarly. 
We introduce the following subcategories of $\dD _{X/Y}$. 
\begin{defi}\label{tec}
For $1\le i \le N$ and $k\in \mathbb{Z}$, we define 
$\aA _{(i,k)}$ by 
$$\aA _{(i,k)}\cneq \left(
\oPPer (X/Y_i) \cap \dD _{X/Y} \right) \otimes \oO _X(kD_i).$$
\end{defi}
The following lemma will be used in the next section.
 We give the proof in Section 6. 
\begin{lem}\label{pers}

(i) For $p=-1, 0$, we have 
$$\Coh _{\cup _{j\neq i}C_j}(X) \cup  \pPPer (\dD _{X/Y_i})
\subset \pPPer (X/Y_i)\cap \dD _{X/Y}, $$
as subcategories of $\dD _{X/Y}$. Moreover
the right hand side coincides with the smallest extension 
closed subcategory of $\dD _{X/Y}$ which contains the left hand side.

(ii) We have the following equality of subcategories of $\dD _{X/Y}$:
$$\iPPer (X/Y_i) \cap \dD _{X/Y} =
\aA _{(i,1)}. $$

(iii) Any object in
 $\aA _{(i,k)}$ is given by a successive extension of 
 the objects in 
$$\Coh _{\cup _{j\neq i}C_j}(X) , \quad \oO _{C_i}(k-2)[1], \quad 
\oO _{C_i}(k-1).$$
Here $C_i$ has a reduced scheme structure. 

\end{lem}
\begin{rmk}\emph{As we see in Section 6, the assumption that
$X_0$ is smooth is needed in the proof of Lemma~\ref{pers} (ii), (iii). 
If otherwise $C_i$ in Lemma~\ref{pers} (iii) may 
not be reduced, so the arguments in the next
section (especially Lemma~\ref{partial}) does not work.}
\end{rmk}

\subsection{Fourier-Mukai transforms and standard equivalences}
The Fourier-Mukai transform is a useful tool in studying 
derived categories. 
Let $g\colon W \to Y$ be another crepant resolution. 
We introduce the relative version of the Fourier-Mukai transform. 
\begin{defi}
We say the equivalence $\Phi \colon D^b (W)\to D^b (X)$ 
is of Fourier-Mukai type over $Y$ if there exists an object 
$\eE \in D^b (X\times W)$, which is supported on 
$X\times _Y W$, such that $\Phi$ is written as 
$$\Phi \cong \Phi _{W\to X}^{\eE} \cneq \dR p_{X\ast}(p_W ^{\ast}(\ast) \dotimes \eE).$$
Here $p_X$, $p_W$ are corresponding projections from 
$X\times W$. The object $\eE \in D^b (X\times W)$ is called the 
kernel of $\Phi$. 
\end{defi}

The following theorem plays an important role in this paper. 

\begin{thm}\emph{\bf{\cite[Theorem 1.1]{Br1}}} \label{standard}, 
\emph{\bf{\cite[Proposition 4.2]{Ch}}}
Let $g\colon W\to Y$ be another crepant small resolution, and 
$\phi \colon W\dashrightarrow X$ be the birational map. Assume 
for a $g$-ample divisor $H$ on $W$,  the divisor
$-\phi _{\ast}H$ is $f$-ample. Then the functor 
$$\Phi _{W \to X}^{\oO _{W \times _{Y} X}} \colon 
D^b (W) \longrightarrow D^b (X),$$
gives an equivalence, and takes $\iPPer (W/Y)$ to 
$\oPPer (X/Y)$. 
\end{thm}
We call the equivalence given in Theorem~\ref{standard}
\textit{standard}.
 
 \begin{rmk}\label{not}
 \emph{ Let us apply Theorem~\ref{standard}
 to the sequence
  $X_i ^{\dag}\to Y_i \leftarrow X$ as in the diagram (1).
  Then there exists a standard equivalence
  $$
  D^b (X_i ^{\dag}) \lr D^b (X),$$
  which takes $\iPPer (X_i ^{\dag}/Y_i)$ to $\oPPer (X/Y_i)$ but 
  not necessary takes $\iPPer (X_i ^{\dag}/Y)$ to $\oPPer (X/Y).$}
  \end{rmk}
 
 We introduce the set $\FM (X)$ as follows. 

\begin{defi}\label{FM}
We define
$\FM (X)$ to be the set of pairs $(W, \Phi)$, where 
$g\colon W\to Y$ is another crepant small resolution and 
$\Phi \colon D^b (W) \to D^b (X)$ is an equivalence of 
derived categories which satisfies the following:
there exists a factorization of the birational map $\phi \colon W
\dashrightarrow X$,
$$W=X^{n}\dashrightarrow X^{n-1}\dashrightarrow \cdots \dashrightarrow 
X^1 \dashrightarrow X^0 =X, $$
and equivalences of Fourier-Mukai type over $Y$, 
$\Phi ^j \colon D^b (X^j) \to D^b (X^{j-1})$ such that 
$\Phi \cong \Phi ^1 \circ \cdots \circ \Phi ^n$. Each 
$\Phi ^j$ is one of the following:
\begin{itemize}
\item \emph{type I :}$X^j = X^{j-1}$ and $\Phi ^j \cong \otimes \lL$ 
for $\lL \in \Pic (X^j)$. 
\item \emph{type I\hspace{-.1em}I :} 
$X^j \dashrightarrow X^{j-1}$ is a flop
at a single rational curve as in the diagram (1),  and 
$\Phi ^j$ is a standard equivalence. 
\end{itemize}
\end{defi} 
Note that for $(W, \Phi) \in \FM(X)$, $\Phi$ restricts to the 
equivalence, 
$$\Phi \colon \dD _{W/Y} \lr \dD _{X/Y},$$
because the kernel of $\Phi$ is supported on $X\times _Y W$. 
Therefore $\Phi$ induces an isomorphism between 
$K(\dD _{W/Y})$ and $K(\dD _{X/Y})$. 
For the description of $K(\dD _{X/Y})$, we have the following lemma. 
We give the proof in Section 6. 

\begin{lem} \label{Kg}
(i) 
The class $[\oO _x] \in K(\dD _{X/Y})$ does not depend on 
a choice of $x\in C$, 
and $K(\dD _{X/Y})$ is described by 
  the direct sum, 
  $$K(\dD _{X/Y})= \mathbb{Z}[\oO _x] \oplus \bigoplus _{i=1}^N 
  \mathbb{Z}[\oO _{C_i}(-1)].$$
  
  (ii) For $(W, \Phi)\in \FM(X)$, $\Phi$ takes the class $[\oO _w] \in 
  K(\dD _{W/Y})$ to $[\oO _x] \in K(\dD _{X/Y})$ for closed 
  points $w\in W$ and $x\in X$. 
  \end{lem}

Finally we give the following lemma, which relates Fourier-Mukai 
transforms to chern characters. The proof will also be given in 
Section 6. 
\begin{lem}\label{chern}
For $(W, \Phi)\in \FM(X)$ and $\lL \in \Pic (W)$, one has 
$$\ch _1 \Phi (\lL) =\phi _{\ast}\ch _1 \lL +\ch _1 \Phi (\oO _W).$$
Here $\phi \colon W\dashrightarrow X$ is the birational map. 
\end{lem}

\section{The chamber structures on the normalized stability conditions}

In this section, we prove Theorem~\ref{cham}.
As in the introduction, let $\Stab (X/Y)$ be the space of stability 
conditions on $\dD _{X/Y}$. 

\subsection{Constructions of stability conditions} 
Here we construct stability conditions, which corresponds to 
the neighborhood of the large volume limit. 
By Proposition~\ref{tstru}, in order to 
give a stability condition, we find
a heart of a bounded t-structure $\aA \subset \dD _{X/Y}$ and a 
stability function $Z\colon K(\dD _{X/Y}) \to \mathbb{C}$. 
As a heart of a bounded t-structure, we take the subcategory
$$\Coh _C(X) =\Coh (X) \cap \dD _{X/Y} \subset \dD _{X/Y}.$$
In order to construct stability functions, 
 let us take an element 
  $\beta + i \omega \in N^1(X/Y)_{\mathbb{C}}$. 
 We define $Z_{(\beta, \omega)}\colon K(\dD _{X/Y})\to \mathbb{C}$ as 
 follows:
 \begin{align*}
 Z_{(\beta, \omega)}(E)&\cneq
  -\int e^{-(\beta +i \omega)}\ch (E)\\
  &= -\ch _3 (E)+(\beta +i \omega)\ch _2 (E).
  \end{align*}
  
  \begin{lem}\label{nor}
  Assume $\beta +i \omega \in A(X/Y)_{\mathbb{C}}$. Then the pair 
  $$\sigma _{(\beta, \omega)}\cneq (Z_{(\beta, \omega)}, \Coh _C (X))$$ 
  determines a point of $\Stab (X/Y)$. 
  \end{lem}
  \textit{Proof}.
  We use Proposition~\ref{tstru} and Remark~\ref{filt}.
   Let us take a non-zero 
  object $E\in \Coh _C (X)$. 
   If $\dim \Supp (E)=1$, then 
  $\Imm Z_{(\beta, \omega)}(E)>0$ since $\omega$ is ample. 
  If $\dim \Supp (E)=0$, then $Z_{(\beta, \omega)}(E)\in \mathbb{R}_{<0}$. 
 Therefore it suffices to show the Harder-Narasimhan property. 
 
 By Remark~\ref{filt}, we check $\Coh _C(X)$ is noetherian. 
 But this is obvious, because a finitely generated module over 
 a noetherian ring is noetherian. Thus it is enough to check 
 that there are no infinite sequence, 
$$ \cdots \subset E_{j+1}\subset E_j \subset \cdots \subset E_2 \subset E_1,$$
for $E_j \in \Coh _C(X)$ 
with $\phi (E_{j+1})>\phi(E_j)$ for all $j$. 
Assume the sequence as above exists. 
For $1\le i \le N$, let $\eta _i \in X$ be the generic point of $C_i$. 
Then since $E_j$ are supported on $C$, one has 
$$0\le \length _{\oO _{X, \eta _i}}E_{j+1}
 \le \length _{\oO _{X, \eta _i}}E_j < \infty.$$
 Therefore we may assume $\length _{\oO _{X, \eta _i}}E_j$ is 
 constant for all $i$ and $j$. Then if we take the exact sequence, 
 $$0 \lr E_{j+1} \to E_j \to G_j \to 0,$$
 then $G_j$ is zero dimensional, thus 
 $Z_{(\beta, \omega)}(G_j)\in \mathbb{R}_{\le 0}$.
 This implies $\phi (E_j)\ge \phi (E_{j+1})$, thus we get a contradiction. 
 $\quad \square$
 
  \hspace{3mm}
  
 The stability conditions
  $\sigma _{(\beta, \omega)}$ are contained in the following 
  subspace:
  $$\Stabn (X/Y)\cneq \{ \sigma =(Z,\pP)\in \Stab (X/Y)
  \mid Z([\oO _x])=-1 \},$$
  for closed points $x\in C$. 
  We call $\Stabn (X/Y)$ the set of 
  \textit{normalized stability conditions}. 
  By Lemma~\ref{Kg} (i),   
 the following map gives an isomorphism of the vector spaces:
 $$N^1 (X/Y)_{\mathbb{C}} \ni \beta +i \omega \longmapsto Z_{(\beta, \omega)}\in
 \{ Z\in \Hom (K(\dD _{X/Y}), \mathbb{C}) \mid Z([\oO _x])=-1 \}.$$
  Hence the map sending stability conditions to their 
  central charges restricts to give the map,
  $$\zZ _{X}\colon \Stabn (X/Y) \lr N^1 (X/Y)_{\mathbb{C}}.$$
We give the following definition. 
\begin{defi}\label{LV}
We define $U_{X}\subset \Stabn (X/Y)$ as follows,
$$U_{X}\cneq \{ \sigma _{(\beta, \omega)}\in \Stabn (X/Y) \mid 
\beta +i \omega \in A (X/Y)_{\mathbb{C}}\}.$$
Also define $$\Stabn ^{\circ}(X/Y) \subset \Stabn (X/Y),$$ 
 to be the connected component 
of $\Stabn (X/Y)$ which 
contains $U_{X}$. 
\end{defi}
Note that $\zZ _{X}$ restricts to a homeomorphism between 
$U_{X}$ and $A (X/Y)_{\mathbb{C}}$.

\subsection{Codimension one boundaries of $U_{X}$}
Here we give the descriptions of codimension one boundaries of 
the domain $U_X$. Let us take a point,
$$\sigma =(Z_{(\beta, \omega)}, \aA) \in \overline{U}_X.$$
Note that by the continuity of $\zZ _X$, we have
$$\beta +i \omega =\zZ _{X}(\sigma) \in \overline{A}(X/Y)_{\mathbb{C}}.$$
We say $\sigma$ \textit{lies in the codimension one boundary}
if $\zZ _X (\sigma)\in \wW _i$ for some $1\le i\le N$. 
The purpose here is to find a heart of a t-structure $\aA$ 
such that the pair $(Z_{(\beta, \omega)}, \aA)$ gives a stability 
condition, and lies in the codimension one boundary. Recall that 
we defined $\wW _{i,k}$ in Definition~\ref{amp} and 
$\aA _{(i,k)}$ in Definition~\ref{tec}.

\begin{lem}\label{partial}
Let us take $\beta +i \omega \in \wW _{i,k}$ for 
some $k\in \mathbb{Z}$.  Then the 
pair 
$$\sigma _{(\beta, \omega)}\cneq (Z_{(\beta, \omega)}, 
\aA _{(i,k)}),$$
gives a stability condition on $\dD _{X/Y}$. 
\end{lem}
\textit{Proof}. 
First we check the following,  
\begin{align}Z_{(\beta, \omega)}(\aA _{(i,k)}) 
\subset \{ \mathbb{R}_{>0} \exp (i\pi \phi) \mid 
0 < \phi \le 1 \}.\end{align}
By
 Lemma~\ref{pers} (iii), $\aA _{(i,k)}$ is the smallest
extension closed subcategory of $\dD _{X/Y}$ which contains
$$\Coh _{\cup _{j\neq i}C_j}(X) , \quad \oO _{C_i}(k-2)[1], \quad 
\oO _{C_i}(k-1).$$
 Hence it is enough to check (2) for the above objects. 
 Let us take $E\in \Coh _{\cup _{j\neq i}C_j}(X)$. If 
 $\dim \Supp (E)=1$, then $\Imm Z(E) >0$ because 
 $\omega \cdot C_j >0$ for $j\neq i$. If $\dim \Supp (E)=0$, 
 one has $Z_{(\beta, \omega)}(E) \in \mathbb{R}_{<0}$. Therefore 
 (2) holds for $E\in \Coh _{\cup _{j\neq i}C_j}(X)$. For 
 other two objects, one has 
 \begin{align*}
 & Z_{(\beta,\omega)}(\oO _{C_i}(k-2)[1])= k-1-\beta \cdot C_i \in 
 \mathbb{R}_{<0},\\
 & Z_{(\beta, \omega)}(\oO _{C_i}(k-1))=-k +\beta \cdot C_i \in 
 \mathbb{R}_{<0}.
 \end{align*} 
 Therefore (2) holds, and it is enough to check 
 $\aA _{(i,k)}$ is noetherian and the pair $(Z_{(\beta, \omega)}, 
 \aA _{(i,k)})$ satisfies
Harder-Narasimhan property, by Remark~\ref{filt}. 
 
First we check $\aA _{(i,k)}$ is noetherian. 
Assume there exists an infinite sequence
$$E_1 \subset E_2 \subset \cdots E_n \subset \cdots \subset E$$
in $\aA _{(i,k)}$. 
We put $E_{\infty}=E$. 
By Lemma~\ref{pers} (iii), 
the complex $E_n$ is concentrated in degree zero at the generic point
of $C_j$ for $j\neq i$. Hence we may assume the generic length of $E_n$
at $C_j$ are constant for all $j\neq i$. Then let us take 
the exact sequence
$$0\to E_{1}\to E_{n} \to G_n \to 0$$
 in $\aA _{(i,k)}$. 
 Then $G_n$ is supported on $C_i$, hence 
 $$G_n \in \oPPer (\dD _{X/Y_i})\otimes 
 \oO _{X}(kD_i).$$ 
 We obtain the sequence 
 $$G_1 \subset G_2 \subset \cdots \subset G_{\infty},$$
 in $\oPPer (\dD _{X/Y_i})\otimes 
 \oO _X (kD_i)$.  
 By Proposition~\ref{simple} $\oPPer (\dD _{X/Y_i})\otimes \oO _X (kD_i)$
  is finite length, hence 
 this sequence terminates. 

 Finally we check the Harder-Narasimhan property. 
Let us assume there exists an infinite sequence 
$$\cdots \subset E_n \subset \cdots \subset E_2 \subset E_1, $$
in $\aA _{(i,k)}$
such that $\phi (E_{n+1})>\phi (E_n)$ for all $n$. 
Then $\Imm Z_{(\beta, \omega)}(E_n)>0$.  
Again we may assume generic length of $E_n$
at $C_j$ are constant for $j\neq i$. If we take  
the exact sequence
$$0\to E_{n+1}\to E_n \to G_n \to 0$$
in $\aA _{(i,k)}$,  
 then $G_n$ is supported on 
$C_i$, hence $Z_{(\beta, \omega)}(G_n)\in \mathbb{R}_{\le 0}$. 
But this implies   
$\phi (E_n)\ge (E_{n+1})$, hence a contradiction. 
$\quad \square$

\hspace{3mm}

The following proposition gives the complete 
description of the codimension one boundaries of $U_{X}$. 
\begin{prop}\label{close}
For $\beta +i \omega \in \wW _i$, there exists a point 
$\sigma \in \overline{U}_X$ with $\zZ _X (\sigma)=\beta +i \omega$
if and only if $\beta +i \omega \in \wW _{i,k}$ for some $k \in \mathbb{Z}$. 
If $\beta +i \omega \in \wW _{i,k}$ for $k\in \mathbb{Z}$, the 
point $$\sigma _{(\beta, \omega)} =(Z_{(\beta, \omega)}, 
\aA _{(i,k)})\in \Stabn (X/Y),$$
in Lemma~\ref{partial} gives the point of $\overline{U}_X$. 
\end{prop}
\textit{Proof}. 
First assume there exists $\sigma \in \overline{U}_X$ 
with $\zZ _X (\sigma)=\beta +i \omega \in \wW _i$. 
It is easy to see that $\oO _{C_{i}}(k-1)$ is stable with respect to 
all $\sigma '\in U_X$ and $k\in \mathbb{Z}$.
 By the comment in~\cite{Brs1} after~\cite[Proposition 8.1]{Brs1}, 
$\oO _{C_{i}}(k-1)$ is at least semistable in $\sigma$. 
This implies 
$$Z_{(\beta, \omega)}(\oO _{C_i}(k-1))=-k +\beta \cdot C_i \neq 0,$$
for all $k\in \mathbb{Z}$. Therefore we have $\beta \cdot C_i \notin 
\mathbb{Z}$. 

Conversely assume that $\beta \cdot C_i \in (k-1, k)$ for 
some $k\in \mathbb{Z}$. 
In this case, we have the stability condition 
$$\sigma \cneq 
\sigma _{(\beta, \omega)}=(Z_{(\beta, \omega)}, \aA _{(i,k)}),$$
constructed in Lemma~\ref{partial}. Therefore it is enough 
to check $\sigma \in \overline{U}_X$. 
We leave it to the next Lemma, and give the proof in Section 6, 
because it requires some technical arguments. $\quad \square$

\begin{lem} \label{diff}
$\sigma _{(\beta, \omega)}=(Z_{(\beta, \omega)}, \aA _{(i,k)}),$
constructed in Lemma~\ref{partial} is contained in $\overline{U}_X$. 
\end{lem}
Now we give the following definition. 
\begin{defi}\label{ik}
We define $\partial{U}_X (i,k)$ for $1\le i \le N$ and $k\in \mathbb{Z}$ 
by 
$$\partial{U}_X (i,k)= \{ (Z_{(\beta, \omega)},\aA_{(i,k)})
\in \overline{U}_X 
\mid \beta +i \omega \in  \wW _{i,k} \}.$$
\end{defi}

\subsection{Other regions by Fourier-Mukai transforms}
Here we construct other regions using $\FM(X)$ defined in
Definition~\ref{FM}. 
Let us take $u=(W, \Phi)\in \FM (X)$. Recall that $\Phi$ gives an 
equivalence $\Phi \colon D_{W/Y} \to D_{X/Y}$. Then 
for any $\sigma =(Z, \pP) \in \Stab (W/Y)$, one gets the 
stability condition $\Phi _{\ast}\sigma =(Z', \pP ')$ by 
$$Z'(E)\cneq Z\circ \Phi ^{-1}(E), \quad 
\pP '(E) \cneq \Phi \pP (\phi), $$
for $E\in \dD _{X/Y}$ and $\phi \in \mathbb{R}$.
Obviously, the map
$$\Phi _{\ast} \colon \Stab (W/Y) \lr \Stab (X/Y)$$
gives a homeomorphism. Furthermore by Lemma~\ref{Kg} (ii),
$\Phi _{\ast}$ restricts to the homeomorphism,
$$\Phi _{\ast} \colon \Stabn (W/Y) \lr \Stabn (X/Y).$$
\begin{defi}\label{other}
For $u=(W, \Phi)\in \FM(X)$, we define the open set $U(W, \Phi)$
in $\Stabn (X/Y)$ to be 
$$U(W, \Phi) \cneq \Phi _{\ast}U_W \subset \Stabn (X/Y),$$
and
 the map $\phi ^{u}\colon N^1 (W/Y)_{\mathbb{C}} \to N^1 (X/Y)_{\mathbb{C}}$
to be 
$$\phi ^u (D)=\ch _1 \Phi (\oO _W)+ \phi _{\ast}D. $$
Here $\phi$ is the birational map $\phi \colon W \dashrightarrow X$. 
\end{defi}
We have introduced the map $\phi ^u$ because 
the following lemma holds.

\begin{lem}\label{com}
For $u=(W, \Phi)\in \FM(X)$, we have the following 
commutative diagram:
$$\xymatrix{
\Stabn (W/Y) \ar[r]^{\Phi _{\ast}}\ar[d]_{\zZ _{W}} & \Stabn (X/Y) \ar[d]^{\zZ _{X}} \\
N^1 (W/Y)_{\mathbb{C}} \ar[r]^{\phi ^{u}}& N^1 (X/Y)_{\mathbb{C}}.
}$$
\end{lem}
\textit{Proof}. Let us take $\sigma =(Z_{(\beta, \omega)}, \pP)
\in \Stabn (W/Y)$, and $\lL \in \Pic (W)_{\mathbb{C}}$ 
with $\ch _1(\lL)=\beta +i \omega$. Then for $F\in \dD _{W/Y}$, 
 the value $Z_{(\beta, \omega)}(F)$ is written as $-\chi (\lL, F)$ by 
Riemann-Roch theorem. Therefore for $E\in \dD _{X/Y}$, one
calculates $Z_{(\beta, \omega)}\Phi ^{-1}(E)$ as follows, 
\begin{align*}
Z_{(\beta, \omega)}\Phi ^{-1}(E) & = -\chi (\lL, \Phi ^{-1}(E)) \\
&= -\chi (\Phi (\lL), E) \\
&= -\ch _0 \Phi (\lL) \ch _3 E + \ch _1 \Phi (\lL)\ch _2 E \\
&= -\ch _3 E + (\phi _{\ast}\ch _1 \lL + \ch _1 \Phi (\oO _W)) \ch _2 E \\
&= -\ch _3 E + \phi ^u (\beta +i \omega)\ch _2 E.
\end{align*}
Here the fourth isomorphism follows from Lemma~\ref{chern}. 
The above equality exactly means $\zZ _X \circ \Phi _{\ast}= \phi ^v 
\circ \zZ _W$. $\quad \square$

\hspace{3mm}

We have the following corollary. 
\begin{cor}
$\zZ _X$ restricts to give a homeomorphism between 
$U(W, \Phi)$ and $\phi _{\ast}A(W/Y)_{\mathbb{C}}$. 
\end{cor}

\subsection{Chambers in $\Stabn (X/Y)$}
Here we investigate the relationship of the regions 
$U(W, \Phi)$. 
\begin{prop}\label{reg}
The regions $U(W, \Phi)$ satisfy the 
following:

(i)
$U(W, \Phi) \cap U(W', \Phi ')\neq \emptyset$ if and only if
there exists $\lL \in \Pic (W)$ such that 
$$W\cong W' \quad \mbox{ and } \quad \Phi ^{-1}\circ \Phi ' \cong \otimes \lL.$$

(ii) $U(W, \Phi) \cap U(W', \Phi ') = \emptyset$ but
$\overline{U}(W, \Phi)\cap \overline{U}(W', \Phi ')\neq 
\emptyset$ in a codimension one wall 
if and only if $W \dashrightarrow W'$ is a flop at a single rational 
curve, 
and we have 
$$\Phi ^{'-1}\circ \Phi  \cong \otimes \lL '\circ \Phi '' \circ \otimes 
\lL $$
with $\Phi ''$ a standard equivalence 
$$\Phi ''\colon D^b (W) \lr D^b (W'), $$
and $\lL \in \Pic (W)$, $\lL '\in 
\Pic (W')$. 
Moreover all points in a codimension one wall of $\overline{U}(W, \Phi)$ 
are contained in $\overline{U}(W', \Phi ')$ with $(W', \Phi')$ as above. 
\end{prop}
\textit{Proof}. 
(i)
We may assume $W'=X$ and $\Phi ' =\id$. 
Suppose $U(W, \Phi)\cap U (X,\id)\neq \emptyset$. 
Then $\Phi$ takes $\Coh _C (W)$ to $\Coh _C (X)$.  Note that simple objects 
of these categories are skyscraper sheaves. Hence for $p\in W$, 
we have $\Phi (\oO _p) \cong \oO _{\phi(p)}$ for some $\phi(p)\in X$. 
Since $\Phi$ is Fourier-Mukai type, $\phi$ gives an isomorphism 
 $W\cong X$ and    
$\Phi \cong \otimes \lL$ for some $\lL \in \Pic (X)$. 
(\textit{cf.}~\cite[Theorem 3.2]{Kawlog}).

(ii)
Again we may assume $W'=X$ and $\Phi ' =\id$. 
First let us consider a flop $\phi _i \colon 
X_i ^{\dag} \dashrightarrow X$ as in the diagram (1). 
Let us consider the standard equivalence,
$$\Phi _i \colon D^b (X_i^{\dag}) \lr D^b (X).$$
 Then recall that
 $\Phi _i$ takes 
$\iPPer (X_i ^{\dag}/Y_i)$ to $\oPPer (X/Y_i)$, by Theorem~\ref{standard}. 
Combined with
 Lemma~\ref{pers} (ii), $\Phi _i$ restricts to give the equivalence:
$$\Phi _i \colon 
\left( \oPPer (X_i ^{\dag}/Y_i) \cap \dD _{X_i ^{\dag}/Y} \right) 
\otimes \oO _{X_i^{\dag}}(D_i^{\dag})
 \longrightarrow \oPPer (X/Y_i) \cap \dD _{X/Y}. $$
Here $D_i ^{\dag}$ is a divisor on $X_i ^{\dag}$ with 
$D_i ^{\dag}\cdot C_i ^{\dag}=1$, and $D_i ^{\dag} \cdot C_j ^{\dag}=0$
for other exceptional curves $C_j ^{\dag}$, $j\neq i$.

By Lemma~\ref{close}, both sides appear as 
t-structures of some points in codimension one boundaries 
of $U_{X_i ^{\dag}}$ and $U_{X}$ respectively.  
On the other hand if $\beta \in N^1 (X_i ^{\dag}/Y)$ satisfies 
$\beta \cdot C_i ^{\dag}\in (0, 1)$, 
then $\phi _{i\ast}\beta \cdot C_i \in (-1, 0)$. 
Therefore Lemma~\ref{close} together with Lemma~\ref{com} 
gives the homeomorphism, 
$$\Phi _{i\ast} \colon \partial{U}_{X^{\dag}_i}(i,1) \lr 
\partial{U}_X (i, 0).$$
Here recall that $\partial{U}_{X}(i,k)$ is defined in Definition~\ref{ik}.
Similarly, we have a homeomorphism, 
 $$(\otimes \oO _X(D_i))_{\ast} \colon 
 \partial{U}_X (i, k) \lr \partial{U}_{X}(i, k+1).$$
 Combining these, for $\lL \in \Pic (X_i ^{\dag})$
  and $\lL ' \in \Pic (X)$
 there exists a homeomorphism, 
 $$(\otimes \lL '\circ \Phi _{i} \circ \otimes \lL )_{\ast} \colon 
 \partial{U}_{X^{\dag}_i}(i,l) \lr 
 \partial{U}_{X}(i,k),$$
 for some $k$ and $l$, and conversely for all $k\in \mathbb{Z}$, there exist 
 $\lL \in \Pic (X_i ^{\dag})$ and $\lL '\in \Pic (X)$ which gives
 the homeomorphism as above. This gives the ``if'' part and the 
 last statement. 
 
 Conversely, assume that $\overline{U}(W, \Phi)$ and $\overline{U}_X$
 intersects at a codimension one boundary, but 
 $U(W, \Phi) \cap U_X =\emptyset$. 
  Then, because $\Stabn (X/Y)$
 is a complex manifold, there exists a pair 
 $(X_i ^{\dag}, \Phi^{\dag})$ with $\Phi ^{\dag}
  = \otimes \lL' \circ \Phi _{i} \circ
  \otimes \lL$ given as above, such that 
  $$U(W, \Phi) \cap U(X_i ^{\dag}, \Phi^{\dag}) \neq \emptyset.$$
  Then (i) implies $W \cong X^{\dag}_i$ and $\Phi ^{\dag}$ is given 
  as $\otimes \lL' \circ \Phi _{i} \circ
  \otimes \lL ''$ for some $\lL '' \in \Pic (W)$. 
 $\quad \square$
 
\subsection{Proof of Theorem~\ref{cham}}
Now we give the proof of Theorem~\ref{cham}. 
\begin{thm}\label{decom}
We have a disjoint union of chambers:
$$\mM _{}\cneq 
\bigcup _{(W, \Phi)\in \FM (X)}
U_{}(W, \Phi) \subset \Stabn ^{\circ}(X/Y),$$
such that two chambers are either equal or disjoint. 
Moreover we have 
$\overline{\mM}_{}=\Stabn ^{\circ}(X/Y)$. 
\end{thm}
\textit{Proof}. 
By the definition, $(W, \Phi)\in \FM(X)$ is obtained as a 
successive sequence of tensoring line bundles and standard equivalences. 
Therefore 
Proposition~\ref{reg} implies $U(W, \Phi)$ are contained in the 
same connected component, i.e. 
$$\mM _{} \subset \Stabn ^{\circ}(X/Y).$$
Also by Proposition~\ref{reg} (i), 
 $U_{}(W, \Phi)\cap U_{}(W', \Phi ')\neq \emptyset$
implies $U_{}(W, \Phi)=U_{}(W', \Phi ')$. 
It remains to show $\overline{\mM}_{}=\Stabn ^{\circ}(X/Y)$. 
Take $\sigma \in \Stabn ^{\circ}(X/Y)$ and 
$\sigma _0 \in U_{X}$. Let 
$$\gamma \colon [0,1]\lr 
\Stabn ^{\circ}(X/Y),$$ be a path such that $\gamma (0)=\sigma _0$ and 
$\gamma (1)=\sigma$. Note that by Theorem~\ref{model}, the chambers 
in $N^1 (X/Y)_{\mathbb{C}}$ given by 
$\phi _{\ast}A (W/Y)_{\mathbb{C}}$ for crepant small resolutions
$g\colon W\to Y$ are finite. 
Because the map $\zZ _{X}$ is a local homeomorphism by Theorem~\ref{lois},
 we may assume 
$\zZ _{X}(\gamma)$ passes through only codimension one 
walls of these chambers
by deforming $\gamma$ if necessary.
Then there exists a sequence, 
$$0< t_1 < t_2 \cdots <t_{n-1}<1,$$
such that $\zZ _X (\gamma(t_k))$ is contained in one of the codimension 
one walls, and $\zZ _X(\gamma ((t_{k-1}, t_k)))$ are contained in 
the chambers. Then $\gamma ((0, t_1))$ is contained in $U_X$ and 
by Proposition~\ref{reg} (ii), 
$\gamma ((t_1, t_2))$ is contained in some $U(W, \Phi)$ with 
$(W, \Phi)\in \FM(X)$. Repeating this argument, we have 
$\gamma (t)\in \overline{\mM} _{}$ for all 
$t\in [0,1]$. In particular, $\sigma =\gamma (1) \in \overline{\mM}$. 
$\quad \square$

\section{Stability conditions and covering spaces}
In this section, we give the proof of Theorem~\ref{cham2}. 
\subsection{Root lattices and Weyl reflections}\label{root}
First we introduce the root lattice associated to the exceptional 
locus of $f\colon X\to Y$. 

\begin{defi}
We define the $\mathbb{Z}$-module $\Lambda _f$ by 
$$\Lambda _f \cneq \bigoplus _{i=1}^N \mathbb{Z}e_i, $$
and introduce the pairing 
$(,) \colon \Lambda _f \times \Lambda _f \to \mathbb{Z}$
by 
$(e_i,e_i) =-2$, $(e_i,e_j) =1$ if $i\neq j$ and $C_i \cap C_j 
\neq \emptyset$, and $(e_i,e_j)=0$ if otherwise. 
\end{defi}
It is well known that the pairing $(,)$ on $\Lambda _f$ is negative 
definite. 
We define the sets $V(\Lambda _f)$, $S(\Lambda _f)$ by 
$$V(\Lambda _f)\cneq \{ v\in \Lambda _f \mid 
(v,v)=-2 \}, \quad S(\Lambda _f)=\{ e_i \mid 1\le i\le N \}.$$
An element $v\in V(\Lambda _f)$ is called a \textit{root}, 
and $e_i \in S(\Lambda _f)$ is called a \textit{simple root}. 
It is well-known that any element $v\in V(\Lambda _f)$ is
written as 
$$v=\sum _{i=1}^N r_i e_i,$$
with $r_i \in \mathbb{Z}_{\ge 0}$ for all $i$, or
 $r_i \in \mathbb{Z}_{\le 0}$ for all $i$. For $v\in V(\Lambda _f)$, 
 we associate the reflection, 
 $$r_v \colon \Lambda _f \ni x \longmapsto x +(v,x)v \in \Lambda _f.$$
 It is well-known that the set $V(\Lambda _f)$ is obtained by applying
 reflections $\{ r_{e_i} \} _{1\le i\le N}$ to the set $S(\Lambda _f)$.   
Note that we have a natural identification of 
vector spaces, 
$$\Lambda _f \otimes _{\mathbb{Z}}\mathbb{R} \cong N_1 (X/Y),$$
by identifying $e_i$ with $C_i$. Then the dual of $r_v$ gives the 
reflection, 
$$r_v ^{\ast} \colon N^1 (X/Y) \lr N^1 (X/Y),$$
and is called the \textit{Weyl reflection}. 
Let us consider a partial resolution $X \stackrel{f'}{\to}Y' \to Y$. 
Then there is a natural inclusion 
$\Lambda _{f'} \hookrightarrow \Lambda _f$,
and it is obvious that
$$V(\Lambda _{f'})=V(\Lambda _f) \cap \Lambda _{f'}.$$
Let $g\colon W\to Y$ be another crepant small resolution. 
Then the dual graphs of 
the exceptional locus of $f$ and $g$ are identified, since
they are identified under simple flops.
(In our case we flops at single rational curves since $Y$ is 
complete, thus this statement is obvious. However 
we note that the explicit constructions of flops~\cite[Theorem 2.4]{Kof} 
would also imply this without the assumption that $Y$ is complete.)
Under the assumption that $X_0$ is smooth, it is possible 
to choose the identification of the dual graphs canonically (does not 
depend on a choice of a decomposition of $W \dashrightarrow X$ into 
flops) because 
the birational map $W\dashrightarrow X$ induces an isomorphism between 
minimal resolutions, $g^{-1}(Y_0) \stackrel{\cong}{\to}X_0$. 
Hence we have an
identification, 
$\psi \colon \Lambda _g \cong \Lambda _f$ which preserves the 
set of simple roots. By taking its dual, we obtain an isomorphism, 
$$\psi ^{\ast} \colon N^1 (X/Y) \lr N^1 (W/Y).$$
By composing with the strict transform, we get the isomorphism, 
$$\phi _{\ast} \circ \psi ^{\ast} \colon N^1 (X/Y) \lr N^1 (W/Y) \lr 
N^1 (X/Y).$$
If $\phi \colon W\dashrightarrow X$ is a flop at $C_i$, the above 
map 
is calculated in~\cite[Theorem 6.3]{Rei}, and it
coincides with $r_{e_i}^{\ast}$. 
Here we note that in~\cite[Theorem 6.3]{Rei}, 
the assumption 
that a
 smooth surface $C\subset X_0$ with $(C)_{X_0}^2 =-2$ exists is  
 needed. 
 Combined with Lemma~\ref{com}, we obtain the following.

\begin{lem}\label{weyl}
We have the map, 
$$\Theta \colon \FM(X) \ni u \longmapsto \phi ^u \circ \psi ^{\ast} \in \Aff
 (N^1 (X/Y)_{\mathbb{C}}).$$
Here the right hand side is the group of affine transformations
on $N^1 (X/Y)_{\mathbb{C}}$. 
Its image coincides with the
 subgroup generated by Weyl reflections 
with respect to simple roots of $\Lambda _f$ and 
translations by the lattice $\Pic (X)\subset N^1 (X/Y)$. 
\end{lem}

\subsection{Proof of Theorem~\ref{cham2}}
Before giving the proof of Theorem~\ref{cham2}, 
we introduce some notations. First note that the following set
$$G\cneq \{ (W, \Phi) \in \FM(X) \mid W=X \},$$
is a subgroup of the group of autoequivalences of $D^b (X)$. 
Then the following map 
$$ G \ni (X, \Phi) \longmapsto \ch _1 \Phi (\oO _X) \in N^1 (X/Y),$$
is a group homomorphism, by the same argument of the proof of 
Lemma~\ref{chern}. We define $G^{\circ}\subset G$
 to be the kernel of the above map. Next for
 $v\in V(\Lambda _f)$, define $H_v \subset N^1 (X/Y)_{\mathbb{C}}$ to be
 $$H_v \cneq \{ \beta +i \omega \in N^1 (X/Y)_{\mathbb{C}} \mid
 (\beta +i \omega) v \in \mathbb{Z} \}.$$
Now we prove Theorem~\ref{cham2}. 
 \begin{thm}\label{gal}
 We have the following map:
 $$\zZ _X\colon 
 \Stabn ^{\circ}(X/Y) \lr N^1 (X/Y)_{\mathbb{C}}\setminus 
 \bigcup _{v\in V(\Lambda _f)}H_v.$$
 It is a regular covering map with Galois group 
 equal to $G^{\circ}$. 
 \end{thm}
 \textit{Proof}. 
 We divide the proof into 3 steps. 
 \begin{step}
 We have 
 $\Imm (\zZ _X)\subset N^1 (X/Y)_{\mathbb{C}}\setminus 
 \bigcup _{v\in V(\Lambda _f)}H_v$. 
 \end{step}
\textit{Proof}.
Note that for $u =(W, \Phi) \in \FM(X)$ with $g\colon W\to Y$, 
$\phi ^u$ restricts to an 
isomorphism, 
$$\phi ^u \colon N^1 (W/Y)_{\mathbb{C}} \setminus 
\bigcup _{v\in V(\Lambda _g)}H_v \lr N^1 (X/Y)_{\mathbb{C}}\setminus 
 \bigcup _{v\in V(\Lambda _f)}H_v.$$
Therefore by Theorem~\ref{decom} and Lemma~\ref{com},
 it suffices to show that for $\sigma =(Z_{(\beta, \omega)},\pP)
\in \overline{U}_{X}$, 
we have $(\beta +i \omega)v \notin \mathbb{Z}$ for all 
$v\in V(\Lambda _f)$. 
First from Proposition~\ref{close}, we have 
$(\beta +i \omega)v\notin \mathbb{Z}$ for simple roots $v$. 
Next suppose
$\omega$ is written as $\omega =f^{'\ast}\omega '$ for 
a partial resolution $X \stackrel{f'}{\to}Y' \to Y$, 
and $\omega '\in A(Y'/Y)$. 
Then $(\beta +i \omega)v$ could be an integer only if 
$$v\in V(\Lambda _f) \cap \Lambda _{f'} =V(\Lambda _{f'}).$$
Take $v\in V(\Lambda _{f'})$. Then by the argument of Lemma~\ref{weyl}, 
there exist flops over $Y'$,
$$\phi \colon W=X^n \dashrightarrow X^{n-1} \dashrightarrow \cdots 
\dashrightarrow 
X^0 =X,$$
and a simple root $v' \in S(\Lambda _{g'})$ such that 
$\phi _{\ast}H_{v'}=H_v$. Here $g'$ is a structure morphism $g'\colon W\to Y'$.
Note that 
$$\omega \in f^{'\ast}A(Y'/Y) \subset \phi _{\ast}\overline{A}(W/Y).$$
Using Theorem~\ref{decom} and the fact that $\zZ _{X}$ is a local 
isomorphism, we can find a region $U(W, \Phi)$
such that 
$\sigma \in \overline{U}(W, \Phi)$. 
Since $\zZ _{W}(\overline{U}_{W})\subset N^1 (W/Y)_{\mathbb{C}}
\setminus H_{v'}$, Proposition~\ref{com}
implies 
$$\beta +i \omega =\zZ _X(\sigma) \in \zZ _{X}(\overline{U}(W, \Phi))\subset 
N^1 (X/Y)_{\mathbb{C}}\setminus H_v.$$
This implies $(\beta +i \omega)v \notin \mathbb{Z}$. 
$\quad \square$

\begin{step}
The map
$\zZ _X$ is surjective. 
 \end{step}
 \textit{Proof}. 
By Lemma~\ref{com}, it suffices to show the surjectivity on 
 $\overline{A}(X/Y)_{\mathbb{C}}\setminus \bigcup _{v\in V(\Lambda _f)}H_v$. 
 Let us take $\beta +i \omega \in \overline{A}(X/Y)_{\mathbb{C}}$ 
 and suppose $\omega =f^{'\ast}\omega '$ for 
 a partial resolution $X \stackrel{f'}{\to}Y' \to Y$, i.e. 
 $\beta +i \omega$ lies in the wall 
 $\wW \cneq N^1 (X/Y)\oplus if^{'\ast}A (Y'/Y)$.
 For a root $v\in V(\Lambda _{f'})$, let 
 $\overline{H}_v$ be 
 $$\overline{H}_v \cneq \{ \beta \in N^1 (X/Y') \mid \beta \cdot v \in \mathbb{Z} \} \subset N^1 (X/Y').$$
 $\overline{H}_v$ is a real codimension one hypersurface and we have the 
 fiber space structure
 $$\Pi \colon 
 \wW \setminus \bigcup _{v\in V(\Lambda _f)}H_v \lr
 N^1 (X/Y')\setminus \bigcup _{v\in V(\Lambda _{f'})}\overline{H}_v 
 ,$$
 with fiber $f^{'\ast}A (Y'/Y)_{\mathbb{C}}$. 
 Let $C(X/Y')$ be 
 one of the connected component of the right hand side, which 
 contains the following set, 
 $$\{ \beta \in N^1 (X/Y') \mid -\varepsilon < \beta \cdot v <0 \mbox{ for all }v\in S(\Lambda _{f'}) \}, $$
 for sufficiently small $0< \varepsilon \ll 1$. 
  We are now going to argue that it suffices to check the surjectivity 
  on $\Pi ^{-1}(C(X/Y'))$.
  
 Let $\gG$ be the group of Weyl reflections generated by the reflections
 by the simple roots 
  $S(\Lambda _{f'})$, acting on $N^1 (X/Y')$. 
 Let $\Pic (X/Y')\subset N^1 (X/Y')$ be the lattice consisting 
 of numerical classes of line bundles. 
 Then the action of 
 $\gG$ preserves $\Pic (X/Y')$. Hence we can consider
  the semi-direct product
 $\gG \ltimes \Pic (X/Y')$ with its action on  
 $N^1 (X/Y')\setminus \bigcup _{v\in V(\Lambda _{f'})}\overline{H}_v$
 given by the following formula:
 $$(g, \lL)\ast D =g(D)+\lL. $$
 Here $g \in \gG$, $\lL \in \Pic (X/Y')$ and $D\in N^1 (X/Y')$.    
 This action also induces the action on the set of connected 
 components of 
 $N^1 (X/Y')\setminus \cup _{v\in V(\Lambda _{f'})}\overline{H}_v$, 
 which can be easily checked to be transitive. 
 
 On the other hand, let us use Lemma~\ref{weyl}. 
 Then for each $(\phi, \lL)\in \gG \ltimes \Pic (X/Y')$, 
 we can find $(W, \Phi)\in \FM (X)$ such that
 \begin{itemize}
 \item A successive sequence 
 $$W=X^{n}\dashrightarrow X^{n-1}\dashrightarrow \cdots \dashrightarrow 
X^1 \dashrightarrow X^0 =X, $$
of Definition~\ref{FM} is defined over $Y'$, i.e. 
the structure morphism $X^j \to Y$ factors as $X^j \to Y' \to Y$. 
  In particular there exists
 a morphism $g' \colon W\to Y'$ and
 $\Theta (W, \Phi)$ preserves the wall $\wW$. 
 \item  We have 
 $\Pi \circ \Theta (W, \Phi)|_{\wW} = (\phi, \lL) \circ \Pi$. 
 \end{itemize}
 Also note that for any lattice 
 isomorphism $\psi \colon \Lambda _{g'} \to 
 \Lambda _{f'}$ which preserve the simple roots,
  the induced morphism $\psi ^{\ast}\colon N^1 (X/Y') \to 
 N^1 (W/Y')$ takes $C(X/Y')$ to $C(W/Y')$. 
 Therefore using Lemma~\ref{com}, we may assume
 $\Pi (\beta +i \omega)$
  lies in $C(X/Y')$. 
 In this region, let us take the pair:
 $$\sigma \cneq (Z_{(\beta, \omega)}, \oPPer (X/Y') \cap \dD _{X/Y}).$$
 By the same proof as in Lemma~\ref{partial} and Proposition~\ref{close}, 
 we can conclude $\sigma \in \Stabn ^{\circ}(X/Y)$. 
 $\quad \square$
 
 \begin{step}
 The map
 $\zZ _X$ is a regular 
 covering map with Galois group $G^{\circ}$. 
 \end{step}
 \textit{Proof}.
 First for $\sigma \in U_X$ and $u=(X, \Phi) \in G^{\circ}$, one has 
 $$\Phi _{\ast}u \in U(X, \Phi) \subset \Stabn ^{\circ}(X/Y).$$
 Therefore the action of $G^{\circ}$ on the set of stability conditions 
 preserve the connected component $\Stabn ^{\circ}(X/Y)$. Moreover 
 by the definition of $G^{\circ}$, one has $\phi ^u =\id$. Therefore 
 Lemma~\ref{com} implies $G^{\circ}$ acts on $\Stabn ^{\circ}(X/Y)$ 
 as deck transformations. This action is also faithful. In fact, 
 assume that for $(X, \Phi)\in G^{\circ}$, 
 one has $\Phi _{\ast}\sigma =\sigma$ for some
  $\sigma \in \Stabn ^{\circ}(X/Y)$. By Theorem~\ref{decom}, we may 
  assume $\sigma \in U(W, \Phi ')$ for some $(W, \Phi ') \in
  \FM (X)$. Then 
  $$U(W, \Phi ')\cap U(W, \Phi \circ \Phi ') \neq \emptyset.$$
  By Lemma~\ref{reg} (i), $\Phi$ is given by tensoring line bundle. 
  Because $\ch _1 \Phi (\oO _X)=0$, we have $\Phi =\id$. 
 
 Next 
 let us take $\sigma, \sigma '\in \Stabn ^{\circ}(X/Y)$ such that 
 $\zZ _X(\sigma)=\zZ _X(\sigma ')$. 
 By Theorem~\ref{decom}, 
 there exist some regions $U(W, \Phi)$, $U(W, \Phi ')$ such that 
 $\sigma \in \overline{U}(W, \Phi)$ and 
 $\sigma ' \in \overline{U}(W, \Phi ')$. 
  Let us choose 
 a sequence $\{\sigma _k \}_{k=1}^{\infty} \subset U(W, \Phi)$ which  
 converges to $\sigma$. Then there exists 
 another sequence $\{ \sigma _k ' \} _{k=1}^{\infty} \subset U(W, \Phi ')$ 
 which converges to $\sigma '$ and
 $\zZ _X(\sigma _k)=\zZ _X(\sigma _k ')$. 
Let us take $\lL \in \Pic (W)$ to be 
$$\lL = \det \Phi ^{'-1} (\oO _X) \otimes \{ \det \Phi ^{-1}(\oO _X)\} ^{-1}.$$
 Then by Lemma~\ref{open} below, 
 the functor 
 $$\Phi ''\cneq \Phi ' \circ \otimes \lL \circ \Phi ^{-1} \colon D^b (X) 
 \to D^b (X),$$
 satisfies $\Phi ''(\oO _X)|_{X^{\circ}}\cong \oO _{X^{\circ}}$, 
 where $X^{\circ}=X\setminus C$. Because $C\subset X$ is codimension 
 two, this implies $\det \Phi ''(\oO _X) \cong \oO _X$, hence 
 $\ch _1 \Phi '' (\oO _X)=0$. 
 Therefore $(X, \Phi '') \in G^{\circ}$ and 
 $$\zZ _X (\Phi ''_{\ast}(\sigma _k))=\zZ _X (\sigma _k)=\zZ _X(\sigma _k').$$
 Note that $\Phi ''_{\ast}$ takes 
 $U(W, \Phi)$ to $U(W, \Phi')$,  
 thus $\Phi ''_{\ast}(\sigma _k)$ is contained in 
 $U(W, \Phi ')$. Since
  $\zZ _X$ is injective on $U(W, \Phi ')$, it follows that 
 $$\Phi ''_{\ast}(\sigma _k)=\sigma _k'. $$
 Taking the limit $k\to \infty$, we get 
 $\Phi _{\ast}''(\sigma)=\sigma '$. 
 This shows $\zZ _X$ is a regular 
 covering 
 map with Galois group $G^{\circ}$. $\quad \square$. 
 
 \hspace{3mm}
 
 Theorem~\ref{gal} provides some 
 information of the group of autoequivalences of $\dD _{X/Y}$. 
 Let $\Aut ^{\circ}(\dD _{X/Y})$ be the group of autoequivalences $\Phi$
 of $\dD _{X/Y}$, such that 
 \begin{itemize}
 \item $\Phi$ preserves the class $[\oO _x] \in K(\dD _{X/Y})$ 
 \item $\Phi$ preserves the connected component $\Stabn ^{\circ}(X/Y)$. 
 \end{itemize}
 Then we have the following. 
 \begin{cor}\label{aut}
 (i)
We have the group isomorphism
$$G =\{ (W, \Phi)\in \FM ^{\circ}(X) \mid 
W=X \} \cong \Aut ^{\circ} (\dD _{X/Y}). $$

(ii)
There exists the following exact sequence:
$$\pi _1 (N^1 (X/Y)_{\mathbb{C}}\setminus 
\cup _{v\in V(\Lambda _{f})}H_v )
\lr \Aut ^{\circ} (\dD _{X/Y}) \lr \Pic (X) \lr 0.$$
Here the last map is defined by taking $\Phi$ to $\det \Phi (\oO _X)$. 
\end{cor}
\textit{Proof}. 
(i) Let us take $\Phi \in \Aut ^{\circ}(\dD _{X/Y})$ and $\sigma \in U_X$. 
Then $\Phi _{\ast}\sigma \in \Stabn ^{\circ}(X/Y)$. By Theorem~\ref{decom}, 
$\Phi _{\ast}\sigma \in \overline{U}(W, \Phi ')$ for some 
$(W, \Phi ')\in \FM (X)$. By deforming $\sigma$ if necessary, we may 
assume $\Phi _{\ast}\sigma \in U(W, \Phi ')$. Then we have 
$$\Phi _{\ast}\sigma \in U(X, \Phi) \cap U(W, \Phi ') \neq \emptyset.$$
Thus by Proposition~\ref{reg} (i),
 we have $W\cong X$ and $\Phi \cong \otimes \lL
\circ 
\Phi '$ for some $\lL \in \Pic (X)$. Therefore 
$(X, \Phi)\in \FM(X)$. 

(ii) This follows from (i) and Theorem~\ref{gal}.
$\quad \square$ 

\begin{rmk}
\emph{If we knew $\Stabn ^{\circ}(X/Y)$ is simply connected, 
we have the exact sequence:
$$0\lr \pi _1 (N^1 (X/Y)_{\mathbb{C}}\setminus 
\cup _{v\in V(\Lambda _{f})}H_v )
\lr \Aut ^{\circ} (\dD _{X/Y}) \lr \Pic (X) \lr 0.$$}
\end{rmk}

\subsection*{Example}
Let $Y=\Spec \mathbb{C}[[x,y,z,w]]/(xy-zw)$ and 
$f\colon X\to Y$ be the blowing up at the ideal $(x,z)$ as in
the introduction. Then Theorem~\ref{gal} describes $\Stabn ^{\circ}(X/Y)$
as a covering map, 
$$\Stabn ^{\circ}(X/Y) \lr \mathbb{C}\setminus \mathbb{Z}.$$
In this case $\Pic (X)=\mathbb{Z}$, and the action of $n\in \Pic (X)$ on 
$N^1 (X/Y)_{\mathbb{C}}=\mathbb{C}$ is given by the translation, 
$\beta +i \omega \mapsto n+\beta +i \omega$. 
Therefore Corollary~\ref{aut} implies, 
\begin{align*}\Stabn ^{\circ}(X/Y) / \Aut ^{\circ}(\dD _{X/Y}) & \cong 
(\mathbb{C}\setminus \mathbb{Z}) / \mathbb{Z} \\
& \cong \mathbb{P}^1 \setminus \{ \mbox{three points} \}.
\end{align*}
Thus we have obtained the same picture as in~\cite[Figure 2]{Asp}.

\begin{rmk}\emph{If we don't assume the existence of a smooth surface
$X_0$,
the description may change. For example, there is an 
example of a crepant small resolution $f\colon X\to Y$
which contracts a
 $(1,-3)$-curve and its scheme theoretic fiber is 
non-reduced~\cite[Lemma 5.16]{Rei}. In this case we have to 
delete extra rational points from
 $\mathbb{C}\setminus \mathbb{Z}$, 
according to the proof of Lemma~\ref{partial}. Thus 
Theorem~\ref{gal} is false in this case.
At this time, the author does not know how to describe $\Stabn ^{\circ}(X/Y)$
in a beautiful way in a general situation involving such cases. 
}
\end{rmk}

\section{Some technical lemmas}
In this section, we give the postponed proofs. 
\subsection{Proof of Lemma~\ref{pers}}
\textit{Proof}. 
(i) For simplicity we show the case of $p=0$. 
The objects of $\oPPer (\dD _{X/Y_i})$ are contained in the right 
hand side by the definition. Take $E\in \Coh _{\cup _{j\neq i}C_j}(X)$. 
Then $\dR g_{i\ast}E \in D^b(Y_i)$
 is a sheaf, hence Lemma~\ref{per} applied to 
$g_i \colon X\to Y_i$ implies 
$E\in \oPPer (X/Y_i)$. Since $E$ is also supported on $C$, $E$ is 
contained in the right hand side.

Conversely let $\bB \subset \dD _{X/Y}$ the smallest extension closed
subcategory which contains the left hand side. 
Let us take $A\in \oPPer (X/Y_i) \cap \dD _{X/Y}$. Then by Lemma~\ref{per}(i),
we have 
the distinguished triangle,
$$
 H^{-1}(A)[1] \lr A \lr H^0 (A) \lr H^{-1}(A)[2].$$
The conditions of Lemma~\ref{per} applied to $g_i \colon X\to Y_i$ shows that
$H^{-1}(A)[1]$ is also contained in $\oPPer (X/Y_i)$. By Lemma~\ref{per}(ii),
one has $g_{i\ast}H^{-1}(A)=0$. 
Therefore $H^{-1}(A)$ must be supported on $C_i$, hence  
$$H^{-1}(A)[1]\in \oPPer (\dD _{X/Y_i}) \subset \bB.$$
Therefore it is enough to check $H^0 (A)\in \bB$. 
 Let $I_{C} \subset I_{C_i}\subset \oO _X$ be the ideal sheaves 
 which define $C$, $C_i$ respectively. Since $H^0(A)$ is 
 supported on $C$, it is annihilated by $I_C ^n$ for some $n \in \mathbb{N}$.
 We have the exact sequence, 
 $$0 \to I_{C_i}^n /I_C ^n \lr \oO _X /I_C ^n \lr \oO _X /I_{C_i}^n \lr 0.$$
 Note that $I_{C_i}^n /I_C ^n$ is supported on $\cup _{j\neq i}C_j$
 and $\oO _X/I_{C_i}^n$ is supported on $C_i$. 
 By tensoring $H^0(A)$, we get the exact sequence
 $$0 \lr A_1 \lr H^0 (A) \lr A_2 \lr 0,$$
 where $A_2 =H^0 (A) \otimes \oO _X /I_{
 C_i}^n$ and 
 $A_1$ is the image of the map $H^0(A)\otimes I_{C_i}^n /I_C ^n \to H^0(A)$.
 Note that $A_1 \in \Coh _{\cup_{j\neq i}C_j}(X) \subset \bB$
  and $A_2 \in \Coh _{C_i}(X)$. 
On the other hand, 
since $R^1 g_{i\ast}H^0 (A)=0$ by Lemma~\ref{per}(ii),
the long exact sequence associated to the
above sequence gives
$R^1 g_{i\ast}A_2=0$. 
Therefore $A_2$ satisfies the conditions of Lemma~\ref{per}, 
thus $A_2 \in \oPPer (\dD _{X/Y_i})\subset \bB$.
Therefore by the exact sequence above, we have 
$H^0 (A)\in \bB$.

(ii) 
Note that the restriction of $g_i \colon X\to Y_i$ to 
$X_0$ is a minimal resolution of $A_1$-singularity. Thus 
the scheme theoretic fiber of $g_i$ is reduced. Hence
by (i) and Proposition~\ref{simple}, the 
left hand side is the smallest extension closed subcategory of $\dD _{X/Y}$ 
which contains
$$\Coh _{\cup_{j\neq i}C_j}(X) \cup \{ \oO _{C_i}, \oO _{C_i}(-1)[1] \}.$$
Similarly the right hand side is the smallest extension closed subcategory 
of $\dD _{X/Y}$ which contains 
\begin{align*}
&\Coh _{\cup_{j\neq i}C_j}(X) \cup \{ \oO _{C_i}(-2)[1], \oO _{C_i}(-1)\}
\otimes \oO _X (D_i) \\
&=\Coh _{\cup_{j\neq i}C_j}(X) \cup \{ \oO _{C_i}(-1)[1], \oO _{C_i} \}.
\end{align*}
Therefore both sides coincide.

(iii)
Since the scheme theoretic fiber of $g_i$ is reduced, 
 this follows from (i) and Proposition~\ref{simple}. 
$\quad \square$

\subsection{Proof of Lemma~\ref{Kg}}
\textit{Proof}. 
(i)
 First we check $[\oO _x]$ does not depend on a choice of
  $x\in C$. Since $[\oO _x]=[\oO _{C_i}]-[\oO _{C_i}(-1)]$ 
  if $x\in C_i$, $[\oO _x]$ does not depend on $x\in C_i$. 
  Moreover since $C$ is connected, one can conclude
  $[\oO _x]$ does not depend on $x\in C$. 
  
  Next 
  one has 
 $$
  \chi (\oO _X, \oO _x)=1, \quad \chi (\oO _X, \oO _{C_i}(-1))=0, \quad
  \chi (\oO _X(D_j), \oO _{C_i}(-1))=-\delta ^{ij}.
  $$
  Here $\delta ^{ij}=1$ if $i=j$ and zero if $i\neq j$. 
  Therefore $\{ [\oO _x], [\oO _{C_i}(-1)] \}_{1\le i\le N}$ are 
  linearly independent. 
  
  Finally, we check the set $\{ [\oO _x], [\oO _{C_i}(-1)] \}_{1\le i\le N}$
  spans $K(\dD _{X/Y})$. 
  For $E\in \dD _{X/Y}$, we can write 
  $[E]=\sum (-1)^i [H^i(E)]$ in $K(\dD _{X/Y})$. So we assume 
   $E\in\Coh _C (X)$. Then one has a filtration 
  $$E_0 \subset E_1\subset \cdots \subset E_{n-1} \subset E_n =E$$
  such that $F_j =E_{j}/E_{j-1}$ is a $\oO _{C_i}$-module for some $i$. 
  Since $[F_j]$ is written as a sum of $[\oO _{x_i}]$ for a $x_i \in C_i$ 
  and $[\oO _{C_i}(-1)]$, $[E]$ is written as a linear 
  combination of $\{ [\oO _x], [\oO _{C_i}(-1)] \}_{1\le i\le N}$.

 (ii) 
 By the definition of $\FM (X)$, $\Phi$ is decomposed into 
 tensoring line bundles and standard equivalences. 
 Cleary tensoring a line bundle preserves the class 
 $[\oO _x]$. Also a standard equivalence 
 preserve the class $[\oO _x]$, because 
 $\Phi$ takes point sheaves to point sheaves 
 away from a flopped curve. $\quad \square$

\subsection{Proof of Lemma~\ref{chern}}
First we introduce some notations, and prepare a lemma. 
For each crepant small resolution $g\colon W\to Y$, 
let $Y^{\circ}\cneq Y\setminus \{0\}$, 
$W^{\circ}\cneq W\setminus \Ex(g)$ and $j_W \colon W^{\circ} 
\hookrightarrow
W$ be the open immersion. 
Here $\Ex (g)$ is the exceptional locus of $g$. 
Let us take $(W, \Phi)\in \FM (X)$. 
Note that the birational map $\phi \colon W\dashrightarrow X$ 
gives an isomorphism $\phi ^{\circ} \colon W^{\circ} \to X^{\circ}$. 
It seems the following lemma is well-known, but we give the proof for 
lack of reference. 
\begin{lem}\label{open}
There exists an equivalence $\Phi ^{\circ}\colon D^b (W^{\circ})
 \to D^b (X^{\circ})$
such that the following diagram commutes,
$$\xymatrix{
D^b (W) \ar[r]^{\Phi} \ar[d]_{j_W ^{\ast}} & D^b (X) \ar[d]^{j_X ^{\ast}} \\
D^b (W^{\circ}) \ar[r]^{\Phi ^{\circ}}& D^b (X^{\circ}).
}$$
More precisely, $\Phi ^{\circ}(\oO _{W^{\circ}})$ is a line bundle, and 
for $E\in D^b (W^{\circ})$ we have 
$$\Phi ^{\circ}(E) \cong \phi _{\ast}^{\circ}E \otimes \Phi ^{\circ}
(\oO _{W^{\circ}}).$$ 
 \end{lem}
 \textit{Proof}. 
 Let $\eE \in D^b (X\times W)$ be the kernel of $\Phi$. 
 Define $\eE ^{\circ}$ to be 
 $$\eE ^{\circ} \cneq \eE |_{X^{\circ}\times W^{\circ}}
 \in D^b (X^{\circ}\times W^{\circ}).$$
 Then let $\Phi ^{\circ}$ be 
 $$\Phi ^{\circ}\cneq 
 \Phi ^{\eE ^{\circ}}_{W^{\circ}\to X^{\circ}} \colon 
 D^b (W^{\circ}) \lr D^b (X^{\circ}).$$
 We check that the above diagram commutes. 
 Let $p_X^{\circ}$, $p_W ^{\circ}$ be the projections 
 from $X^{\circ}\times W^{\circ}$ onto 
 $X^{\circ}$ and $W^{\circ}$, and
  $p_X ' \colon X^{\circ}\times W \to X^{\circ}$ be the projection. 
  For $E\in D^b (W)$, one has 
  \begin{align*}
  j_X ^{\ast}\Phi (E) & \cong 
  j_X ^{\ast}\dR p_{X\ast}(p_W ^{\ast}(E) \dotimes \eE) \\
  & \cong 
  \dR p_{X\ast}^{'}\{(p_W ^{\ast}(E) \dotimes \eE)|_{X^{\circ}\times 
  W}\}  \\
  & \cong \dR p_{X\ast}^{\circ}
  \{(p_W ^{\ast}(E) \dotimes \eE)|_{X^{\circ}\times W^{\circ}}\} 
  \\
  &\cong \dR p_{X\ast}^{\circ}
  \{ p_{W}^{\circ \ast}j_W ^{\ast}E \dotimes \eE ^{\circ}\} \\
  &\cong \Phi ^{\circ}\circ j_W ^{\ast}(E).
  \end{align*}
  Here the second isomorphism follows from the 
  base change formula 
  for the diagram below, 
  $$\xymatrix{
  X^{\circ} \times W \ar[r] \ar[d]_{p_X ^{'}}
  & X\times W \ar[d]^{p_X} \\
  X^{\circ} \ar[r]^{j_X} & X,\\
  }$$
 and the third isomorphism follows because $\eE$ is supported on 
 $X\times _Y W$. 
 
 Finally we check $\Phi ^{\circ}$ gives an equivalence. 
If $\Phi$ is given by tensoring $\lL \in \Pic (X)$, then $\Phi ^{\circ}$ 
is also given by tensoring $\lL |_{X^{\circ}}$. 
If $\Phi$ is given by the standard equivalence, $\Phi ^{\circ}$ is 
given by $\phi ^{\circ}_{\ast}$.  Therefore in 
these cases, $\Phi ^{\circ}$ is an equivalence. Also if we take 
$(W', \Phi ')\in \FM (W)$, then it is easily checked that 
$$\Phi ^{\circ}\circ \Phi ^{'\circ}\cong (\Phi \circ \Phi ')^{\circ}
\colon D^b (W^{'\circ}) \lr D^b (W^{\circ}).$$
Therefore by the definition of $\FM (X)$, the functor $\Phi ^{\circ}$ 
gives an equivalence. Moreover by the description of $\Phi ^{\circ}$ for 
type I and I\hspace{-.1em}I in Definition~\ref{FM},
 for any $(W, \Phi)\in \FM(X)$, 
$\Phi ^{\circ}$ is given by 
$$\Phi ^{\circ}(E) \cong \phi _{\ast}^{\circ}(E) \otimes \lL, $$
for some $\lL \in \Pic (X^{\circ})$. Applying $E=\oO _{W^{\circ}}$, 
one has $\lL =\Phi ^{\circ}(\oO _{W^{\circ}})$. 
 $\quad \square$

\hspace{3mm}

 \textit{Proof of Lemma~\ref{chern}}.
 
For $E\in D^b (X)$, let $\det (E) \in \Pic (X)$ be the determinant 
line bundle. Because $\ch _1 E=\ch _1 \det (E)$
and $\ch _0 \Phi (\oO _W)=1$, 
it is sufficient to prove
$$\det \Phi (\lL) \cong  \det (\phi _{\ast}\lL \otimes \Phi (\oO _W)).$$
Since $C\subset X$ has codimension two, one has 
$j_{X\ast}j_{X}^{\ast}E \cong E$ for a vector bundle $E$ on $X$. 
Hence it is enough to show the above formula on $X^{\circ}$. 
On the other hand, by Lemma~\ref{open} we have 
\begin{align*}
\Phi (\lL)|_{X^{\circ}} & \cong \Phi ^{\circ} (\lL |_{W^{\circ}}) \\
& \cong \phi _{\ast}^{\circ}(\lL |_{W^{\circ}}) \otimes 
\Phi ^{\circ}(\oO _{W^{\circ}}) \\
& \cong (\phi _{\ast}\lL \otimes \Phi (\oO _W))|_{X^{\circ}}.
\end{align*}
Note that taking determinants and the restrictions to $X^{\circ}$ commute. 
Thus one has 
\begin{align*}
\det (\Phi (\lL))|_{X^{\circ}} & \cong \det (\Phi (\lL)|_{X^{\circ}}) \\
& \cong \det ((\phi _{\ast}\lL \otimes \Phi (\oO _W))|_{X^{\circ}}) \\
& \cong \det (\phi _{\ast}\lL \otimes \Phi (\oO _W))|_{X^{\circ}},
\end{align*}
 and the lemma follows. $\quad \square$

\subsection{Proof of Lemma~\ref{diff}}
Before proving Lemma~\ref{diff}, we prepare some lemmas. 
First we give the following technical lemma on the 
perverse t-structures. 
Recall that by the definition, $\aA _{(i,0)}$ is given by 
$$\aA _{(i,0)}=\oPPer (X/Y_i)\cap \dD _{X/Y}.$$

\begin{lem}\label{ti}
(i) For any object $E\in \aA _{(i,0)}$, 
the following sequence is an exact sequence in 
$\aA _{(i,0)}$. 
$$0\lr H^{-1}(E) [1] \lr E \lr H^0 (E) \lr 0,$$
with $H^{-1}(E)[1] \in \oPPer (\dD _{X/Y_i}).$

(ii) For any object $E\in \Coh _C (X)$, there exists an exact sequence 
$$0 \lr F \lr E \lr G \lr 0,$$
in $\Coh _C(X)$ such that 
$$F\in \aA _{(i,0)} \cap \Coh _C(X), \quad G[1] \in \oPPer (\dD _{X/Y_i}).$$ 
\end{lem}
\textit{Proof}.
(i) By the criterion of Lemma~\ref{per}
applied to $g_i \colon X \to Y_i$, for $E\in \aA _{(i,0)}$
the objects
 $$H^{-1}(E)[1], \quad H^0 (E)$$
  are also contained in $\aA _{(i,0)}$. Because $\aA _{(i,0)}$
   is a heart of a bounded t-structure on $\dD _{X/Y}$, 
the distinguished triangle 
$$H^{-1}(E) [1] \lr E \lr H^0 (E)$$
gives a short exact sequence in $\aA _{(i,0)}$. 
 By Lemma~\ref{per} (ii), we have $g_{i\ast}H^{-1}(E)=0$. 
Therefore $H^{-1}(E)$ is supported on $C_i$. 

(ii) For this proof, we refer to~\cite[Section 3]{MVB}. 
 In fact, according to~\cite[Section 3]{MVB}, 
the category $\oPPer (X/Y_i)$ is obtained from 
$\Coh(X)$ as a \textit{tilting}~\cite{HRS} for a certain
 torsion pair $(\tT _0, \fF _0)$ on $\Coh (X)$,
  with $\tT _0$ a torsion part 
 and $\fF _0$ a free part. 
According to~\cite[Section 3]{MVB}, an object  
$F\in \fF _0$ 
must satisfy $g_{i\ast}F=0$. 
Therefore $F$ is supported on $C_i$, and this shows the torsion 
theory $(\tT _0, \fF _0)$ induce the torsion pair 
$(\tT _0 ', \fF _0 ')$ on 
$\Coh _C(X)$. Clearly corresponding tilting is 
$\aA _{(i,0)}.$ Therefore we can decompose 
$E\in \Coh _C(X)$ into the exact sequence 
$$0 \lr T \lr E \lr F \lr 0,$$ 
such that $T\in \tT _0 '$ and $F \in \fF _0 '$. By the construction of the 
tilting~\cite[Proposition 2.1]{HRS}, $T$ and $F[1]$ are objects in $
\aA _{(i,0)}$. Because $F$ must be supported on 
$C_i$, one has $F[1] \in \oPPer (\dD _{X/Y_i})$. 
$\quad \square$

\hspace{3mm}

We use the assumption that 
$X_0$ is smooth in the following lemma. 
\begin{lem}\label{ch}
There exists a constant $K>0$ such that 
if $E\in \Coh _C(X)$ satisfies $\Hom (E, E)=\mathbb{C}$, 
we have 
$$\ch _2 (E) =\sum _{i=1}^N r_i C_i,$$
with $0\le r_i \le K$. 
\end{lem}
\textit{Proof}. 
By~\cite[Lemma 3.8]{UI}, 
$E$ is a $\oO _{f^{-1}(0)}$-module where 
$f^{-1}(0)$ is a scheme theoretic fiber of $f$
at the closed point $0 \in Y$. 
In particular for a hyperplane section 
$0\in Y_0 \subset Y$ and its pull back $X_0 =f^{-1}(Y_0)$, 
$E$ is an $\oO _{X_0}$-module. 
Let us regard $E$ as an object of $\Coh (X_0)$. 
Note that we have 
$$\Hom _X (E, E)=\Hom _{X_0}(E, E)=\mathbb{C}.$$
Because we assume $X_0$ is smooth, 
we can apply Riemann-Roch theorem on $X_0$, and get 
\begin{align*}
-\ch _1 (E)^2 &= \sum (-1)^i \dim \Ext _{X_0}^{i}(E, E)\\
&\le 2. \end{align*}
Here we have taken the chern character $\ch _1(E)$ in $N_1 (X_0/Y_0)$. 
It is well-known that the pairing 
$$N_1 (X_0/Y_0) \times N_1 (X_0/Y_0) \ni (D_1, D_2) \longmapsto 
D_1 \cdot D_2 \in \mathbb{R},$$
is negative definite. Therefore there are only finite number of possibilities 
for $r_i \in \mathbb{Z}_{\ge 0}$ which satisfy 
$$-\left( \sum _{i=1}^{N} r_i C_i \right) ^2 \le 2,$$
and this shows the lemma. $\quad \square$
 
\hspace{3mm}

Let us take $K>0$ as in Lemma~\ref{ch}, 
 and define $\sS \subset \Coh _C(X)$ to be 
 $$\sS \cneq \{ E \in \Coh _C(X) \mid \ch _2 (E)=\sum _{i=1}^{N}
r_i C_i, \mbox{ with } 0\le r_i \le K \}.$$

\textit{Proof of Lemma~\ref{diff}}.

For $n\in \mathbb{N}$, let $\omega _n$ be
$$\omega _n \cneq \omega +\frac{1}{n}D_i \in A(X/Y).$$
Because $\omega _n$ is ample, one has the stability condition
$$\sigma _n \cneq (Z_{(\beta, \omega _n)}, \Coh _C (X)) \in U_X, $$
constructed in Lemma~\ref{nor}. It suffices to show 
$\sigma _n$ converges to $\sigma$.

For $E\in \aA _{(i,k)}$, we 
denote by $\phi (E) \in (0,1]$ the phase with respect to 
the stability function $Z_{(\beta, \omega)}$ of Proposition~\ref{tstru}. 
Also for a $\sigma$-semistable object $E\in \dD _{X/Y}$, we use 
the same notation for its phase $\phi (E) \in \mathbb{R}$. 
Of course these two notions are compatible, i.e. if 
$E\in \aA _{(i,k)}$ is $\sigma$-semistable, the above 
two phases $\phi (E)$ coincide. Similarly, we use the notation 
$\phi _n(E')$ for the phase of $E' \in \Coh _C(X)$ 
with respect to $Z_{(\beta, \omega _n)}$ 
and $\sigma _n$-semistable object $E' \in \dD _{X/Y}$. 

Note that $\omega _n$ converges to $\omega$. Hence by the 
topology on $\Stab (X/Y)$ introduced in~\cite[Section 6]{Brs1}, and 
using~\cite[Lemma 6.1]{Brs1}, it suffices to show the following:
for any $\varepsilon >0$, there exists $M>0$ such that $n>M$ implies, 
\begin{align}
& \mbox{for any }\sigma \mbox{-stable object }E\in \aA _{(i,k)},
\tag*{} 
\mbox{ if } F_n, F_n' \in \dD _{X/Y}  \mbox{ are the semistable factors}
\\ & \mbox{of }
E \mbox{ in } 
\sigma _n  \mbox{ whose phases are the largest and the smallest respectively},
 \mbox{ one has } \\
&  \lvert \phi _n (F_n)- \phi (E) \rvert < \varepsilon, \qquad 
 \lvert \phi _n (F_n')- \phi (E) \rvert < \varepsilon. \tag*{}
 \end{align}

 For simplicity we treat the case of $k=0$. Other cases are similarly 
discussed. 
 We divide the proof into two steps. 
 
 \begin{sstep}
 There exists $M>0$ such that (3) holds for any 
  $n>M$ and any $\sigma$-stable object $E\in \aA _{(\beta, \omega)}$ 
 with $\phi (E)=1$.
 \end{sstep}
 \textit{Proof}. 
 By Lemma~\ref{pers} (i), 
 any object $E\in \aA _{(i,0)}$ is given by 
 a successive extension of the objects in 
 $$\Coh _{\cup _{j\neq i}C_j} (X), \quad \oO _{C_i}(-2)[1], \quad 
 \oO _{C_i}(-1).$$
 Note that any sheaf $F\in \Coh _{\cup _{j\neq i}C_j}
 X$ 
 with $\dim \Supp (F)=1$ 
 satisfies 
 $\Imm Z_{(\beta, \omega)}(F)>0$ because $\omega \cdot C_j >0$ for 
 $j\neq i$. Therefore if $E\in \aA _{(i,0)}$ 
 satisfies $\Imm Z_{(\beta, \omega)}(E)=0$, 
 $E$ is given by a successive extension of the following objects, 
 $$\{ \oO _x \mid x \in \cup _{j\neq i}C_j \}, \quad
 \oO _{C_i}(-2)[1], \quad \oO _{C_i}(-1).$$
 Therefore any $\sigma$-stable object $E$ of $\phi (E)=1$ is 
 one of the above objects. Note that these objects are also stable in 
 $\sigma _n$. 
 One can easily calculate $\phi _n (\oO _x)=\phi (\oO _x)=1$ and 
 $$
 \lim _{n\to \infty}\phi _n (\oO _{C_i}(-2)[1]) 
 =\phi (\oO _{C_i}(-2)[1])=1, \quad 
 \lim _{n\to \infty}\phi _n (\oO _{C_i}(-1)) =\phi (\oO _{C_i}(-1))=1.$$
 Therefore it is enough to choose $M$ so that 
 $n>M$ implies 
 $$\lvert \phi _n (\oO _{C_i}(-2)[1]) -
 \phi (\oO _{C_i}(-2)[1]) \rvert <\varepsilon, 
 \quad 
 \lvert \phi _n (\oO _{C_i}(-1)) -\phi (\oO _{C_i}(-1)) 
 \rvert <\varepsilon .\quad \square$$

 \begin{sstep}
 There exists $M'>0$ such that (3) holds for any 
  $n>M'$ and any $\sigma$-stable object $E\in \aA _{(i,0)}$ 
 with $0<\phi (E)<1$.
 \end{sstep}
\textit{Proof}. 
Let us take a $\sigma$-stable object
$E\in \aA _{(i,0)}$ with $0<\phi (E)<1$. 
Then we have the exact sequence in $\aA _{(i,0)}$ by 
Lemma~\ref{ti} (i), 
$$0 \lr H^{-1}(E)[1] \lr E \lr H^0 (E) \lr 0.$$
By Lemma~\ref{ti} (i), $H^{-1}(E)$ is supported on $C_i$. Therefore 
$\phi (H^{-1}(E)[1])=1$ and the stability of $E$ implies 
$H^{-1}(E)[1]=0$. Therefore we have 
$$E\in \aA _{(i,0)} \cap \Coh _C(X).$$
Let $F_n \in \Coh _C(X)$ be the semistable factor of $E$ in $\sigma _n$
of the largest phase. We have the exact sequence in $\Coh _C(X)$, 
\begin{align}0 \lr F_n \lr E \lr G_n \lr 0, \end{align}
with 
\begin{align}
\phi _n (F_n) \ge \phi _n (E). 
\end{align}
By Lemma~\ref{ti} (ii), we have the exact sequence in $\Coh _C(X)$,
\begin{align}0 \lr H_n \lr F_n \lr I_n \lr 0,\end{align}
such that 
$$H_n \in \aA _{(i,0)}\cap \Coh _C(X),
 \quad 
I_n [1] \in \oPPer (\dD _{X/Y_i}).$$
Here we assume $H_n \neq 0$ and $I_n \neq 0$. In the case of
 one of $H_n$ or $I_n$ is 
zero, one can argue similarly, and in fact it is easier. 
So we leave it to the reader. 
Combining the exact sequences (4) and (6), one gets the exact sequence, 
\begin{align}0 \lr H_n \lr E \lr G_n' \lr 0.\end{align}
Note that $E\in \aA _{(i,0)}\subset \oPPer (X/Y_i)$. 
Thus Lemma~\ref{per} (2) implies $R^1 g_{i\ast}E=0$. 
Therefore the long exact sequence associated to (7) implies 
$R^1 g_{i\ast}G_n'=0$. Again by Lemma~\ref{per}, we have 
$G_n '\in \oPPer (X/Y_i)$. This means the exact sequence 
(7) is also an exact sequence in $\aA _{(i,0)}$. 
By the stability of $E$, one gets, 
\begin{align}
\phi (H_n) < \phi (E).
\end{align}
Then we apply Lemma~\ref{bound} below.
Note that because $E$ is $\sigma$-stable, we have 
$\Hom (E,E)=\mathbb{C}$. 
Therefore Lemma~\ref{ch} implies $E\in \sS$, and because $H_n$ is 
a subsheaf of $E$, we have 
$H_n \in \sS$. 
Furthermore, since $E$ and $H_n$ are objects of $\aA _{(i,0)}$, 
one has 
$$Z_{(\beta, \omega)}(E)\neq 0, \quad Z_{(\beta, \omega)}(H_n) \neq 0.$$
 Let us take $\varepsilon '=
\frac{1}{3}\varepsilon$. Then by Lemma~\ref{bound}, 
we can find $M'>0$ which does not depend on $E$, such that $n>M'$ implies 
\begin{align}
\lvert \phi _n (E)-\phi (E) \rvert <\varepsilon ', \quad
\lvert \phi _n (H_n)-\phi (H_n) \rvert <\varepsilon '.\end{align} 
Because $I_n [1]\in \oPPer (\dD _{X/Y_i})$, we have 
$$Z_{(\beta, \omega)}(I_n[1]) \in \mathbb{R}<0.$$
Thus we apply Lemma~\ref{bound} 
to get $0< \phi _n (I_n) <\varepsilon '$ for $n>M'$. 
Combined with the sequence (6), we have either 
\begin{align}
0 < \phi _n (F_n) < \varepsilon ' \quad \mbox{ or } \\
\phi _n (H_n) \ge \phi _n (F_n).\end{align}
In the first case, we have 
\begin{align*}
\lvert \phi _n (F_n)-\phi (E) \rvert & 
\le ( \phi _n (F_n)-\phi _n (E) )
 + \lvert \phi _n (E)-\phi (E) \rvert \\
& < 2\varepsilon ' < \varepsilon.
\end{align*}
 Here we have used (5),(9) and (10). In the second case, we have 
 \begin{align*}
 \lvert \phi _n (F_n)-\phi (E) \rvert &  \le ( \phi _n (F_n)-\phi _n (E) )
+ \lvert \phi _n (E)-\phi (E) \rvert  \\
& < (\phi _n (H_n)-\phi _n(E)) +\varepsilon ' \quad \mbox{ from }(11), (9)\\
& < (\phi (H_n)-\phi(E)) +2\varepsilon  ' +\varepsilon ' 
\quad \mbox{ from }(9)\\
& < 3\varepsilon ' =\varepsilon \quad \mbox{ from }(8).
\end{align*}
By the same argument, one can show $\lvert \phi _n (F_n')-\phi (E)\rvert <
\varepsilon$ for the semistable factor $F_n'$ in $\sigma '$
with the smallest phase, when $n>M'$. Therefore the condition (3) holds. 
$\quad \square$

\hspace{3mm}

Here we have used the following lemma:
\begin{lem}\label{bound}
We fix $\beta +i \omega \in \wW _i$ as
 in Proposition~\ref{close}. Then for 
each $\varepsilon '>0$, there exists $M'>0$ such that 
$n>M'$ implies 
$$\sup \left\{ \left\lvert \frac{Z_{(\beta, \omega _n)}(F)}
{Z_{(\beta, \omega)}(F)}-1 \right\rvert : F\in 
\sS , Z_{(\beta, \omega)}(F)\neq 0
\right\} <\varepsilon '.$$
\end{lem}
\textit{Proof}.
One can calculate as follows:
\begin{align}
\left\lvert \frac{Z_{(\beta, \omega _n)}(F)}
{Z_{(\beta, \omega)}(F)}-1 \right\rvert & = 
\frac{1}{n}\cdot  \frac{
\ch _2 F \cdot D_i } { \lvert -\ch _3 F + \beta \ch _2 F +i \omega \ch _2 F
\rvert }.
\end{align}
If $Z_{(\beta, \omega)}(F)\neq 0$, then 
$\omega \cdot \ch _2 (F) \neq 0$ or $-\ch _3 F +\beta \ch _2 F \neq 0$. 
In the former case, one has 
$$(12) \le \frac{1}{n}\cdot \frac{D_i \cdot \ch _2 F}{\omega \cdot \ch _2 F}
,$$
and there are only finite number of possibilities for the values 
$$\left\{ \frac{D_i \cdot \ch _2 F}{\omega \cdot \ch _2 F}
: F\in \sS, \ch _2 F \cdot \omega \neq 0 \right\}.$$
In the latter case, one has 
$$(12) \le \frac{1}{n}\cdot \frac{D_i \cdot \ch _2 F}{\lvert -\ch _3 F +\beta \ch _2 F \rvert}
,$$
and because there are finite number of possibilities for 
the values $\{ \beta \cdot \ch _2 F : F\in \sS \}$, we have 
$$\inf \{\lvert -\ch _3 F +\beta \ch _2 F \rvert : F\in \sS, 
 -\ch _3 F +\beta \ch _2 F \neq 0 \} >0.$$
Therefore there exists a constant $K'$ independent of $F\in \sS$ with 
$Z_{(\beta, \omega)}(F)\neq 0$ such that 
$(12) \le K'/n$. Hence for each $\varepsilon '$, one can take 
$M'$ to be $K'/\varepsilon$. $\quad \square$

Yukinobu Toda, Graduate School of Mathematical Sciences, University of Tokyo

\textit{E-mail address}:toda@ms.u-tokyo.ac.jp

\end{document}